\documentclass[number,citesort,MSNbibl,seceqn,dvips]{arxbj}
\usepackage{accents}

\aid{0}
\volume{18}
\issue{1}
\pubyear{2012}
\firstpage{360}
\lastpage{389}
\doi{10.3150/10-BEJ333}

\makeatletter
\newcommand{\fracb}[2]{(#1)/#2}
\newcommand{\fracc}[2]{#1/(#2)}
\newcommand{\fracd}[2]{(#1)/(#2)}
\newcommand{\cal}{\mathcal}
\newtheorem{theo}{Theorem}[section]
\newtheorem{prop}{Proposition}[section]
\newtheorem{cor}{Corollary}[section]
\newremark{exam}{Example}[section]
\newtheorem{lem}{Lemma}[section]
\newremark{rem}{Remark}[section]
\newproclaim{defin}{Definition}[section]
\renewcommand{\mathring}{\accentset{\circ}}
\newcommand{\card}{\operatorname{card}}
\newcommand{\supp}{\operatorname{supp}}
\newcommand{\Conv}{\operatorname{Conv}}
\newproclaim{NA}{Notation (additional)}
\newproclaim{comment}{Comment}
\newremark{remarks}{Remarks}
\newremark{remark}{Remark}
\makeatother

\begin{document}
\begin{frontmatter}

\title{Asymptotics of the maximal radius of an $L^r$-optimal sequence
of quantizers}
\runtitle{Maximal radius of quantizers}

\begin{aug}
\author{\fnms{Gilles} \snm{Pag\`{e}s}\corref{}\thanksref{e1}\ead[label=e1,mark]{gilles.pages@upmc.fr}} \and
\author{\fnms{Abass} \snm{Sagna}\thanksref{e2}\ead[label=e2,mark]{abass.sagna@upmc.fr}}
\runauthor{G. Pag\`{e}s and A. Sagna}
\address{Laboratoire de Probabilit\'{e}s et Mod\`{e}les al\'{e}atoires,
UMR~7599,
Universit\'{e} Paris 6, case 188, 4, pl. Jussieu, F-75252 Paris Cedex
5, France.
\printead{e1,e2}}
\end{aug}

% HISTORY:
\received{\smonth{7} \syear{2010}}

% ABSTRACT
%
\begin{abstract}
Let $P$ be a probability distribution on $\mathbb{R}^d$ (equipped with
an Euclidean norm
$\vert\cdot\vert$). Let $ r> 0 $ and let $(\alpha_n)_{n \geq1}$ be an
(asymptotically) $L^r(P)$-optimal sequence of $n$-quantizers.
We investigate the asymptotic behavior of the maximal radius sequence
induced by the sequence $(\alpha_n)_{n \geq1}$ defined for every
$n \geq1$ by $\rho(\alpha_n) = \max\{\vert a \vert, a \in\alpha_n \}$.
When $\card(\supp(P))$ is infinite, the maximal radius sequence
goes to
$\sup\{ \vert x \vert, x \in\operatorname{supp}(P) \}$ as $n$ goes to infinity.
We then give the exact rate of convergence for two classes of
distributions with
unbounded support: distributions with hyper-exponential tails and
distributions with polynomial tails.
In the one-dimensional setting, a sharp rate and constant are provided
for distributions with
hyper-exponential tails.
\end{abstract}

% KEYWORDS
%
\begin{keyword}
\kwd{distribution tail}
\kwd{function with regular variation}
\kwd{maximal radius of a quantizer}
\kwd{optimal quantization}
\kwd{Zador theorem}
\end{keyword}

\end{frontmatter}

%s1 ###
\section{Introduction}\label{sec1}

The aim of this paper (which is a part of the second author's Ph.D.
thesis \cite{Sag1}) is to provide some precise upper and lower
bounds for the radius of a sequence of quantizers of an $\mathbb
{R}^d$-valued random vector. Our motivation is that it is a first
attempt toward the elucidation of the geometric structure of an optimal
quantizer in higher dimension.

Quantization has become an important field of information theory since
the early $1940$'s. Nowadays, it plays an important role in digital
signal processing (DSP), the basis of many areas of technology, from
mobile phones to modems and multimedia PCs. In DSP, vector quantization
is the process of approximating a continuous range of values or a very
large set of discrete values by a relatively small set of discrete
values. A common use of quantization is the conversion of a continuous
signal into a digital signal. This is performed in analog-to-digital
converters with a given quantization level.

Recently, optimal vector quantization has become a promising tool in
numerical probability: it is an efficient method to produce grids
optimally fitted to the distribution of a random vector~$X$. This leads
to some cubature formulas that may approximate either expectations (see
\cite{Pag}) or, more significantly, conditional expectations (see
\cite{PagPhaPri1}). This ability to approximate conditional
expectations is the key property called upon in the quantization-based
numerical schemes used to solve some problems arising in finance,
including optimal stopping problems (pricing and hedging American-style
options, see \cite{BalPag,BalPagPri}), the pricing of swing options
(see \cite{BarBouPag,BarBouPag1}), stochastic control problems (see
\cite{CorPhaRun,PagPhaPri}) for portfolio management and nonlinear
filtering (see \cite{PagPha,PhaRunSel}). Other applications, like
some new schemes for the discretization of Zakai and McKean--Vlasov
equations, have also been investigated (see \cite{GobPagPhaPri}).

At this stage, we need to recall some basic facts on optimal
quantization. At this level of generality, we just assume that $\mathbb
{R}^d$ is endowed with a norm $|\cdot |$, possibly not Euclidean.

Let $X \in L^r(\Omega,\mathcal{A},\mathbb{P})$ be an $\mathbb
{R}^d$-valued random vector with distribution $P=\mathbb{P}_X$. The
$L^r(P)$-optimal quantization problem at level $n$ for $X$ consists in
finding the best approximation of $X$ by $q(X)$ for the $L^r(\mathbb
{P})$-norm, where $q$ is a Borel function taking at most $n$ values.
This leads to the following minimization problem:
\[
 \inf\{ \Vert X - q(X) \Vert_r, q\dvtx \mathbb{R}^d \stackrel{\mathrm
{Borel}}{\longrightarrow} \mathbb{R}^d, \card(q(\mathbb{R}^d))
\leq n \},
\]
where $\operatorname{card}(\alpha)$ stands for the cardinality of $\alpha$. The
solution, $e_{n,r}(X),$ of the previous problem is called the
$L^{r}$-optimal mean quantization error induced by $X$ (at level $n$).
Note that, in fact, $e_{n,r}(X)$ only depends on the distribution of
$X$ so that we will occasionally use the notation $e_{n,r}(P)$.
However, for every Borel function $q\dvtx \mathbb{R}^d \rightarrow\alpha$,
\mbox{$\alpha\subset\mathbb{R}^d, \operatorname{card} (\alpha) \leq n$}, we have
\[
\vert X -q(X) \vert\geq d(X,\alpha) := \mathop{\min}_{a\in\alpha} \vert X - a
\vert \qquad\mathbb{P} \mbox{-a.s.}
\]
Consider $\alpha\subset\mathbb{R}^d$ with $\operatorname{card} (\alpha) \leq
n$ (called an $n$-quantizer). Let $(C_a(\alpha))_{a\in\alpha}$ be a
Voronoi partition of $\mathbb{R}^d$ (with respect to the norm $\vert
\cdot\vert$), that is, a Borel partition of $\mathbb{R}^d$ satisfying
for every $a \in\alpha$,
\[
C_a(\alpha) \subset\Bigl\{ x \in\mathbb{R}^d\dvt \vert x-a \vert= \min_{b
\in\alpha} \vert x-b \vert\Bigr\}
\]
and let $ \widehat{X}^{\alpha} = \sum_{a\in\alpha} a \mathbf{1}_{\{
X \in C_a(\alpha)\}}$. Then $ \widehat{X}^{\alpha}$ is a projection on
$\alpha$ following the nearest neighbor rule and satisfying $\vert X -
\widehat{X}^{\alpha} \vert= d(X,\alpha)$ so that one also has
%
%e1.1 ###
\begin{eqnarray} \label{er.quant}
e_{n,r}(X) &=&  \mathop{\mathop{\inf}_{{\alpha\subset\mathbb{R}^d }}}_
{\operatorname{card}(\alpha) \leq n}  \biggl(\int_{\mathbb{R}^d} d(x,\alpha)^r P(\mathrm{d}x)
 \biggr)^{1/r}\nonumber
 \\[-8pt]
 \\[-8pt]
   &=&\inf \{  (\mathbb{E} \vert X -
\widehat{X}^{\alpha} \vert^r  )^{1/r}, \alpha\subset\mathbb{R}^d,
\operatorname{card}(\alpha) \leq n  \}.
\nonumber
\end{eqnarray}
For every $n \geq1$, the infimum in (\ref{er.quant}) holds as a
(finite) minimum {attained by (at least) one so-called \textit{$L^r(P)$-optimal $n$-quantizer} $\alpha^{\star}$ (see, e.g., \cite
{Pag}, Proposition 11 or \cite{GraLus}, Theorem 4.1), also called,
especially when dealing with numerical applications, \textit{the optimal
$n$-grid}. A sequence of $n$-quantizers $ (\alpha_n)_{n \geq1} $ is
$L^r(P)$-optimal if, for every $n \geq1$, $\alpha_n$ is
$L^r(P)$-optimal. A sequence $ (\alpha_n)_{n \geq1} $ is \textit{asymptotically $L^r(P)$-optimal} if
\[
 \int_{\mathbb{R}^d}d(x,\alpha_n)^r P(\mathrm{d}x) = e_{n,r}^r(X) +
\mathrm{o}(e_{n,r}^r(X)) \qquad \mbox{as }  n \rightarrow \infty
\]
($f(x) = \mathrm{o}(g(x))$, as $x \rightarrow\infty$, if $f(x)=\epsilon(x)
g(x)$ with $\lim_{x \rightarrow\infty} \epsilon(x) =0$ for two
$\mathbb R$-valued functions $f$ and $g$).} Moreover, the
$L^r(P)$-optimal mean quantization error $e_{n,r}(X)$ decreases to $0$
as $n$ goes to infinity. As soon as $X$ has a finite $r'$-moment for
some $r' >r$, its rate of convergence to $0$ is ruled by the so-called
Zador theorem.
\begin{theo}[(Zador theorem, see \cite{GraLus,BucWis,ZAD})] Let $X  \in
L^{r'}(\mathbb{P})$ for an $r'>0$, with distribution $P=f \lambda_d+
P_s$ (where $P_s$ denotes the singular part of $P$ with respect to
$\lambda_d$). Then,
%
%e1.2 ###
\begin{equation}\label{Zadorrate}
\forall  r \in(0,r') \qquad\lim_{n} n^{r/d} (e_{n,r}(P))^r = Q_r(P),
\end{equation}
 where
%e1.3 ###
\begin{equation}\label{QrP}
   Q_r(P) = J_{r,d}  \biggl( \int_{\mathbb
{R}^d} f^{\fracc{d}{d+r}}\,\mathrm{d}\lambda_d  \biggr)^{\fracb{d+r}{d}} = J_{r,d}
\Vert f \Vert_{\fracc{d}{d+r}}  \in[0,+\infty),
\end{equation}
 with
 \[
     J_{r,d} = \inf_{n \geq1} n^{r/d}
e_{n,r}^r(U([0,1]^d)) \in(0,+\infty)
\]
and $U([0,1]^d)$ stands for the uniform distribution on $[0,1]^d$.
\end{theo}

Note that $\mathbb{E}\vert X \vert^{r'} < + \infty$ implies $\Vert f
\Vert_{\fracc{d}{d+r}} <+\infty$ and that $J_{r,d} $ depends upon the
norm $|\cdot |$ on~$\mathbb{R}^d$.

Let us come back to our topic of interest, that is, the asymptotic
behavior of the radii of a sequence $(\alpha_n)_{n \geq1}$ of
$L^r$-optimal quantizers. The maximal radius (or simply radius) $\rho
(\alpha)$ of a quantizer $\alpha\subset\mathbb{R}^{d}$ is defined by
\[
\rho(\alpha) = \max\{ \vert a \vert, a \in\alpha \}.
\]
In a one-dimensional setting ($d=1$), one can define the \textit{one-sided} (right) radius of $\alpha$ by removing absolute values in
the above definition. The one-sided left radius is defined as the
opposite of the right radius of $-\alpha$ viewed as a quantizer of $-X$.

From now on, $\vert\cdot\vert$ will denote an Euclidean norm on
$\mathbb{R}^d$, except where explicitly stated otherwise. Except in ambiguous
cases, we will denote $(\rho_n)_{n \geq1}$ for the sequence $(\rho
(\alpha_n))_{n \geq1}$ of radii of $(\alpha_n)_{n \geq1}$.

We will first show that, if the support of $P$, denoted $\operatorname{supp}(P)$, is unbounded, then
$ \lim_{n \rightarrow+\infty}
\rho_n  =  +\infty$ (when $d = 1$, the sequence of one-sided right
radii goes to infinity as soon as $\sup\operatorname{supp}(P) = +\infty$).
The key inequalities to get the upper and lower estimates of the
maximal radius sequence are provided in Theorems~\ref
{prop_princip_limsup} and~\ref{prop_gen_liminf}. In these theorems, we
point out the close connection between the asymptotics of $\rho_n$ and
the generalized survival function of $X$ defined on $\mathbb{R}_+  :=
[0,+\infty)$ by $\bar{F}_r(\xi)  =  \mathbb{E}(\vert X \vert^r
 \mathbf{1} _{\{ \vert X \vert>\xi\}})$. The regular variation index will
play an important role since we elucidate the asymptotic behaviour of
$\rho_n$ (or $\log\rho_n$) from the asymptotic behavior of the
function $-\log\bar{F}_r$ as a regularly varying function. We present
below two typical results obtained for important families of
(essentially radial) distributions: a sharp rate for $\log\rho_n$ for
distributions with polynomial tails and an exact rate for $\rho_n$ for
distributions with hyper-exponential tails (also made sharp when $d =
1$ and $r \ge 1$).

\begin{theo}\label{1.2} Let $P= f \lambda_d$.
\begin{longlist}[(b)]
\item[(a)] \textup{Polynomial tail.} If there exists $K>0$, $\beta
\in\mathbb{R}$, $c> r+d$ and a real number $A>0$ such that
\[
\forall  x \in\mathbb{R}^d \qquad|x|\ge A  \quad \Longrightarrow \quad  f(x) =K
\frac{({\log}\vert x \vert)^{\beta}}{\vert x \vert^c},
\]
then
\[
\lim_{n} \frac{\log\rho_n}{\log n} = \frac{1}{c-r-d} \frac{r+d}{d}.
\]
\item[(b)] \textup{Hyper-exponential tail.} If there exists $K>0$,
$\kappa$, $\vartheta>0$, $c \in\mathbb{R}$ and a real number $A>0$
such that
\[
\forall  x \in\mathbb{R}^d \qquad |x|\ge A  \quad \Longrightarrow \quad  f(x) =K
\vert x \vert^{c} \mathrm{e}^{- \vartheta \vert x \vert^{\kappa}},
\]
 then
\[
   \vartheta^{-1/ \kappa} \biggl (1+\frac{r}{d}
\biggr)^{1/\kappa} \leq \liminf_{n} \frac{\rho_n}{ ( \log n)^{1/\kappa}}
\leq \limsup_{n} \frac{\rho_n}{ (\log n)^{1/\kappa}} \leq2 \vartheta
^{-1/ \kappa} \biggl (1+\frac{r}{d}  \biggr)^{1/\kappa}.
\]
Furthermore, if $d=1$ and $r \geq1$,
\[
\lim_{n} \frac{\rho_n}{ (\log n)^{1/\kappa}}= \vartheta^{-1/ \kappa}
 ( 1+r  )^{1/\kappa}.
\]
\item[(c)] If $f$ has a one-sided polynomial or hyper-exponential
tail, say on $\mathbb{R}_+$, then the maximal radius sequence satisfies
the above asymptotic bounds.
\end{longlist}
\end{theo}

\begin{remarks*}

$\bullet$ Note that the Euclidean norm
appearing in the statement of the above theorem needs to be the one
used to define the radius and the distance between the random vector
and the quantizer. If $X$ has a ${\cal N}(0,I_d)$ distribution, this
norm is the canonical one. As concerns the ${\cal N}(0,\Sigma)$
distribution, the ``reference'' Euclidean norm is $|\cdot |_{\Sigma
^{-1}}$ induced by the inverse $\Sigma^{-1}$ ($\vert x \vert^2_{\Sigma
^{-1}} := x' \Sigma^{-1} x$ for a (column) vector $x \in\mathbb R^{d}$
with $x'$ standing for the transpose of $x$). To derive asymptotic
bounds from such results for the radius \textit{measured in the canonical
Euclidean norm} one needs to use the strong equivalence of the norms,
namely $\frac{1}{\lambda_{\Sigma,\max}}|\cdot |\le|\cdot |_{\Sigma
^{-1}}\le\frac{1}{\lambda_{\Sigma, \min}}|\cdot |$, where $\lambda
_{\Sigma,\max}$ and $\lambda_{\Sigma,\min}$ are the maximum and the
minimum eigenvalues of $\Sigma$, respectively.

$\bullet$ Note that as concerns asymptotic lower estimates,
we propose in Section~\ref{loweriid} an alternative approach based on
random quantization.
\end{remarks*}

The paper is organized as follows. We first give, as a preliminary
result, the limit of the maximal radius for distributions supported by
a set of infinite cardinality. Section~\ref{sec2} is devoted to the upper
estimate of the maximal radius based on the asymptotic estimates of
survival functions of~$X$. Section~\ref{sec3} is devoted to the lower limit
where our results are obtained by two different methods -- one still
based on survival functions and one based on mean random quantization.
In both cases, we strongly rely on recent results obtained in \cite
{GraLusPag} about the $L^s$-behaviour of $L^r$-optimal quantizers when
$r<s<r+d$.

\begin{NA*}   For every $r\ge0$, we
define $ L^{r+}(\mathbb{P}) = {\bigcup}_{\varepsilon>0}L^{r+\varepsilon
}(\mathbb{P})$ and the generalized $r$-survival function $\bar{F}_r(\xi
) = \mathbb{E}  (\vert X \vert^r  \mathbf{1} _{\{ \vert X \vert>
\xi\}}  )$ of a random vector $X \in L^r(\mathbb{P})$. Note that
$\bar{F}_r$ is defined on $\mathbb{R}_{+}$ and takes values in
$[0,\mathbb{E} \vert X \vert^r]$. $\bar{F_0}$ is the regular survival
function denoted $\bar{F}$.

Let $A\subset\mathbb{R}^d$. $\overline{A}$ will stand for its closure,
$\partial A$ for its boundary, $\operatorname{Conv}(A)$ for its convex hull,
$\mathring{A}$ for its interior and $A^{c}$ for its complement.
$[x]$ will denote the integer part of an $x \in\mathbb{R}$.
$B(x,r)$, $r>0$, will denote the open ball with center $x \in\mathbb
{R}^d$ and radius $r\ge0$ and~$d(x,A)$ the distance of $x$ to the set
$A\subset\mathbb{R}^d$. For $x,y \in\mathbb R^{d}$, $(x|y)$ will
denote the inner product of $x$ and $y$ with respect to the specified
Euclidean norm and for two real-valued functions $f$ and $g$, $f(x)
=\mathrm{O}(g(x))$ as $x \to\infty$ if there is a positive real constant $C$
such that $ |f(x)| \leq C |g(x)| $, for all large enough $x$.
\end{NA*}

%-------------------- ASYMPTOTIC OF THE maximal radii
%-----------------------------
%s2 ###
\section{A first preliminary result}\label{sec2}
As a first necessary step we
elucidate the connections between the asymptotics of the maximal radius
sequence and the ``supremum'' of the support of the distribution $P$.
\begin{prop} \label{proplimitRM}
\textup{(a)} Let $\vert\cdot\vert$ be an arbitrary norm on $\mathbb{R}^d$
and $X \in L^r(\mathbb{P})$. Let $(\alpha_n)_{n\geq1}$ be a sequence
of $n$-quantizers such that
$\int_{\mathbb{R}^d}d(x,\alpha_n)^rP(\mathrm{d}x)\to0$ as $n \rightarrow
+\infty$. Then,
%
%e2.1 ###
\begin{equation} \label{hyp6.1}
\liminf_{n} \rho_n \geq \sup\{ \vert x \vert, x \in\operatorname{supp}(P) \}.
\end{equation}

\textup{(b)} Suppose that $\vert\cdot\vert$ is an Euclidean norm
on $\mathbb{R}^d$. If $\card(\supp(P)) = +\infty$, then for any
$L^r(P)$-optimal sequence of $n$-quantizers $(\alpha_n)_{n \geq1}$
%
%e2.2 ###
\begin{equation}
\lim_{n} \rho_n = \sup_{n \geq1} \rho_n = \sup\{ \vert x \vert, x
\in\operatorname{supp}(P) \}.
\end{equation}
\end{prop}

\begin{pf}
(a)
Let $ x \in\supp(P)$ and let $\varepsilon>0$. For every $n
\geq1$,
\begin{eqnarray*}
\Vert d(X,\alpha_{n}) \Vert_r
& \geq& \Vert d(X,B(0,\rho_{n})) \Vert_r   \qquad  \bigl(\mbox{since } \alpha_{n} \subset B(0, \rho_{n})\bigr) \\
& \geq& \bigl\Vert d(X,B(0,\rho_{n}))  \mathbf{1} _{\{X \in
B(x,\varepsilon)\}} \bigr\Vert_r \\
& \geq& d(B(x,\varepsilon),B(0,\rho_{n})) \mathbb{P}\bigl(X \in
B(x,\varepsilon)\bigr)^{1/r}.
\end{eqnarray*}
Consequently, $d(B(x,2\varepsilon),B(0,\rho_{n}))=0$ for large enough
$n$ since $\Vert d(X,\alpha_{n}) \Vert_r \rightarrow0 $ so that
$|x|-2\varepsilon\le\rho_n$, which eventually implies $ \liminf_{n}
\rho_n \geq\vert x \vert$.

(b) We will show first that if $\alpha$ is an $L^r$-optimal quantizer
at level $n$ and if $\card(\supp(P)) \geq n,$ then
%
%e2.3 ###
\begin{equation} \label{Eqlimitesup}
\alpha\subset\overline{\Conv(\supp(P))} \quad \mbox{and}\quad
\rho_n \leq \sup\{ \vert x \vert, x \in\operatorname{supp}(P) \}.
\end{equation}
Note first that if $\alpha$ is $L^r$-optimal at level $n$, then
$\operatorname
{card}(\alpha) = n$ since $\card(\supp(P)) \geq n$ (see \cite
{Pag}, Proposition 11 or \cite{GraLus}, Theorem 4.1). Now, suppose
that there exists $ a \in\alpha\cap (\overline{\Conv(\supp(P))}  )^{c}$ and set
\[
\alpha' = (\alpha\backslash\{ a \}) \cup\{ \Pi(a) \},
\]
where $\Pi$ denotes the projection on the non-empty closed convex set
$\overline{\Conv(\supp(P))}$. The projection is $1$-Lipschitz (see,
e.g., \cite{HirLem}, Chapter III, page~116) and $X$ is $\mathbb{P}\mbox{-a.s. }
  \operatorname{supp}(P)$-valued, hence
%
%e2.4 ###
\begin{equation} \label{EqIneqDistance}
d(X,a) \geq d(\Pi(X),\Pi(a)) \stackrel{\mathbb{P}\mathrm{\mbox{-}a.s.}}{=}
d(X,\Pi(a)).
\end{equation}
It follows that
\[
d(X,\alpha) \geq d(X,\alpha') \qquad\mathbb{P}\mbox{-a.s.}
\]
Since $\alpha$ is $L^r(P)$-optimal at level $n$ and $\operatorname{card}(\alpha
') \leq\operatorname{card}(\alpha)=n$,
\[
\mathbb{E}(d(X,\alpha')^r) = \mathbb{E}(d(X,\alpha)^r)
\]
so that the following two statements hold:
\begin{itemize}
\item $d(X,\alpha') = d(X,\alpha)$  $\mathbb{P}\mbox{-a.s.}$
\item $\Pi(a) \notin\alpha\backslash\{ a \} $ since $\alpha'$ is
$L^r(P)$-optimal (which implies that $\operatorname{card} (\alpha') = n$).
\end{itemize}
On the other hand, it follows from equation (\ref{EqIneqDistance}) that
\[
 \bigl(a-\Pi(a) \vert X-\Pi(a)  \bigr) \leq0 \qquad\mathbb{P} \mbox{-a.s.}
\]
Consequently
\begin{eqnarray*}
\vert X-a \vert^2 - \vert X - \Pi(a) \vert^2 & = & 2\bigl(\Pi(a)-a \vert
X-\Pi(a)\bigr) + \vert a - \Pi(a) \vert^2 \\
& \geq& \vert a-\Pi(a) \vert^2 >0 \qquad\mathbb{P}\mbox{-a.s.}
\end{eqnarray*}
since $a \notin\overline{\Conv(\supp(P))}$. As a consequence
\[
d(X,\alpha') < d(X,\alpha) \qquad\mathbb{P}\mbox{-a.s. on } \bigl\{ X \in
\mathring{C}_{\Pi(a)}(\alpha') \bigr\},
\]
where $ \mathring{C}_{\Pi(a)}(\alpha') = \{\xi \in\mathbb{R}^d, d(\xi
,\Pi(a))< d(\xi,\alpha\backslash\{ a \}) \} $ since the norm is Euclidean.

This implies that $\mathbb{P}(X \in\mathring{C}_{\Pi(a)}(\alpha')) =
0$; if so, $\alpha'\setminus\{\Pi(a)\}=\alpha\setminus\{a\}$ would
clearly be optimal at level $n$ (since $d(X,\alpha)=d(X,\alpha\setminus
\{a\})$ a.s.) with a cardinality equal to $n-1,$ which is impossible
since $e_{n,r}(X)$ decreases (strictly) to $0$ (see again \cite
{GraLus,Pag}). Hence $ \alpha\subset\overline{\Conv(\supp(P))} $.

Now, let us prove that $ \rho_n \leq \sup\{ \vert x \vert, x \in
\operatorname{supp}(P) \}.$ Note first that this assertion is obvious if $\operatorname{supp}(P)$ is unbounded. Otherwise, if $\operatorname{supp}(P)$ is bounded, then
it is compact and so is $\Conv(\supp(P))$. Let $x_0 \in \Conv(\supp(P))$ be such that $\vert x_0 \vert= \sup\{ \vert x \vert,
x \in\Conv(\supp(P)) \}.$ Thus
\[
x_0 = \lambda_0 \xi_1 + (1-\lambda_0) \xi_2, \qquad\xi_1,\xi_2  \in
\operatorname{supp}(P)
\]
and $\lambda\mapsto\vert\lambda\xi_1 + (1-\lambda) \xi_2 \vert$ is
convex so that it attains its maximum at $\lambda = 0$ or $\lambda
= 1$. Consequently $x_0 \in\operatorname{supp}(P)$. Hence $\rho_n  \le
\sup\{ |x |,   x \in\operatorname{supp}(P) \} $, which, combined with~(\ref
{proplimitRM}), yields the conclusion.
\end{pf}

\begin{remark*} Note that (b) follows from the fact that if
$\alpha$ is an $L^r$-optimal quantizer at level $n$, then
%
%e2.5 ###
\begin{equation} \label{AssertRestNorm}
\alpha\subset\overline{\Conv(\supp(P))}
\end{equation}
as soon as $\card(\supp(P)) \geq n$. But this result holds true
only for Euclidean norms on~$\mathbb{R}^d$. For an arbitrary norm, this
assertion may fail. A counterexample is given with the $l_{\infty
}$-norm in \cite{GraLus}, page~25.

Before dealing with the general case we give two examples of
distributions (exponential and Pareto) for which the sharp convergence
rate of the maximal radius sequence can be easily derived from
semi-closed forms established in \cite{ForPag} for their $L^r$-optimal
quantizers.
\end{remark*}

$\rhd$ \textit{Exponential distribution.} Let $r > 0$ and let $P$ be an
exponential distribution with parameter $\lambda  > 0$.
Then
%
%e2.6 ###
\begin{equation} \label{asymp_exp}
\rho_n   =   \frac{r+1}{\lambda} \log n+ \frac{C_{r}}{\lambda
} + \mathrm{O}\biggl (\frac{1}{n}  \biggr),
\end{equation}
where $C_{r}$ is a real constant depending only on $r$.

$\rhd$ \textit{Pareto distribution.} Let $r>0$ and let $P$ be a Pareto
distribution with index $\gamma>r$. Then,
%
%e2.7 ###
\begin{equation} \label{asymp_Par}
\rho_n = K_r n^{\fracd{r+1}{\gamma-r}}  \biggl(1 + \mathrm{O}  \biggl(\frac
{1}{n}  \biggr)  \biggr),
\end{equation}
where $K_r$ is a positive real constant depending only on $r$.

 A short proof of these results is given in the \hyperref[appm]{Appendix}.
These rates will be useful to validate the asymptotic rates obtained by
other approaches.

%--------------------- CONVERGENCE RATE ---------------------------

% ---------------------------------- UPPER ESTIMATE
%----------------------------------
%s3 ###
\section{Asymptotic upper bounds for the radius}\label{sec3}
We investigate in this section the upper rate of convergence of $(\rho
_n)$ to infinity. We next give some definitions and some hypotheses
that will be useful later on.

Let $(\alpha_n)_{n \geq1}$ be an $L^r(P)$-optimal sequence of
quantizers at level $n$. For every $n \geq1$, we denote by $M(\alpha
_n)$ the set of points in $\alpha_{n}$ for which the maximal norm is
reached, namely,
\[
M(\alpha_n) = \Bigl\{a \in\alpha_n \mbox{ such that } \vert a \vert=
{\max_{b \in\alpha_n}} \vert b \vert \Bigr\}.
\]
We will need the following (light) assumption on the distribution $P$:
\[
\mathbf{(H)} \equiv \exists x_0 \in\mathbb{R}^d,
\exists \varepsilon_0>0, \exists  r_0>0  \mbox{ such that }
P(\mathrm{d}x) \geq\varepsilon_0  \mathbf{1} _{B(x_0,r_0)} (x) \lambda_d(\mathrm{d}x),
\]
which means that $P$ is locally lower bounded as a measure by the
Lebesgue measure on a ball. This assumption holds as soon as $P$ has a
density $f$, bounded away from $0$ on a~non-empty open set.

In order to get a sharp estimate for $\rho_n$ for one-dimensional
distributions with hyper-exponential tails, we will need the following
more technical assumption (for $r \in[1,+\infty)$):

$\mathbf{(G_r)} \equiv P=f \cdot\lambda_1,$ where $f>0$ is
non-increasing to $0$ on $[A,+\infty)$, non-decreasing from~$0$ on
$(-\infty,-A]$ for some real constant $A\ge0$ and
%
%e3.1 ###
\begin{equation} \label{Assump_exact_limit}
\lim_{|y| \rightarrow+\infty} \int_1^{+\infty} (u-1)^{r-1} \frac
{f(uy)}{f(y)}\,\mathrm{d}u =0.
\end{equation}
Such an assumption is clearly satisfied by distributions with
hyper-exponential tails, that is, of the form
$ f(x) = K\vert x \vert^{c} \mathrm{e}^{- \vartheta \vert x \vert^{\kappa}}$, $
|x| >A>0$, $\vartheta,  \kappa>0$, $c \in\mathbb{R}$.
Indeed, such a density $f$ is non-increasing outside a compact
interval and we have
\[
 \int_1^{+\infty} (u-1)^{r-1} \frac{f(uy)}{f(y)}\,\mathrm{d}u = y^{-c}\int
_1^{+\infty} (u-1)^{r-1}u^c \mathrm{e}^{-\vartheta y^{\kappa}(u^{\kappa}-1)}\,\mathrm{d}u
\stackrel{y \rightarrow+\infty}{\longrightarrow} 0
\]
by the Lebesgue convergence theorem. A \textit{one-sided version} of
condition $\mathbf{(G_r)}$ can be stated by restricting $f$ on
$[A,+\infty)$ or $(-\infty,-A]$ for some $A\ge0$.

%s3.1 ###
\subsection{Main results on asymptotic upper bounds}\label{sec3.1}
The main result of this section, stated below, makes the connection
between the asymptotic behaviour of $\rho_n$ and that of its survival
function (through some asymptotic ``semi-inverse'' of $-\log\bar F_r$
or $-\log\bar F_r(\mathrm{e}^{\centerdot})$), where $\bar{F}_r(\xi) = \mathbb{E}
 (\vert X \vert^r  \mathbf{1} _{\{ \vert X \vert> \xi\}}  )$
denotes the generalized survival function.

First we need to briefly recall some background on inverse function and
regular variations.

It is clear that the function $\bar{F}_r$ is non-increasing and goes
to $0$ as $\xi\rightarrow+\infty$ (provided $\mathbb{E} \vert X \vert
^{r}<+\infty$). Consequently, $\xi\mapsto-\log\bar{F}_r(\xi)$ is
monotone non-decreasing and goes to ${+}\infty$ as $\xi$ goes to ${+}\infty$.

It is well known that if a function $f$ defined on $(0,+\infty)$ is
non-decreasing to ${+}\infty$, its generalized inverse function
$f^{\leftarrow}$ defined for every $y>0$ by
%
%e3.2 ###
\begin{equation}
f^{\leftarrow}(y)=\inf\{\xi>0, f(\xi) \geq y \}
\end{equation}
is non-decreasing to ${+} \infty$. If, furthermore (see \cite
{BinGolTeu}, Theorem 1.5.12.), $f$ is regularly varying (at ${+}\infty$)
with index\vadjust{\goodbreak} $1/ \delta$, $\delta>0$ (i.e., for every $t>0$, $ \frac
{f(t\xi)}{f(\xi)} \rightarrow t^{1/\delta}$ as $\xi\rightarrow+\infty
$),\vspace*{2pt} then there exists a function $\psi$, regularly varying with index
$\delta$ and satisfying
%
%e3.3 ###
\begin{equation} \label{Asymp_equiv}
\lim_{\xi\to+ \infty} \frac{\psi(f(\xi))}{\xi} = \lim_{y\to+\infty
}\frac{f(\psi(y))}{y} =1 .
\end{equation}
Such a function $\psi$ is called an \textit{asymptotic inverse} of $f$.
It is neither necessarily increasing nor continuous. Moreover, $\psi$
is unique up to asymptotic equivalence at ${+}\infty$ and $f^{\leftarrow
}$ is one version of $\psi$. (By asymptotic equivalence (at ${+}\infty$),
we mean $f\sim g$ if $ \lim_{x \rightarrow+\infty} \frac{f(x)}{g(x)}=1$.)

We show in the theorem below how to derive from the regularly varying
property of a function $\psi_r$ with upper bounds $(-\log\bar
{F}_r)^{\leftarrow}$ or $(-\log\bar{F}_r(\mathrm{e}^{\centerdot}))^{\leftarrow
}$ an asymptotic upper estimate for $\rho_n$ or $\log(\rho_n)$.

\begin{theo} \label{thm_princip_limsup}
Let $r>0$ and let $X \in L^{r}(\mathbb{P})$ with distribution $P$
having an unbounded support and satisfying $(\mathbf{H})$. Let $(\alpha
_n)_{n \geq1}$ be an $L^r(P)$-optimal sequence of $n$-quantizers.
\begin{longlist}[(b)]
\item[(a)] If $\psi_r$ is a non-decreasing function, regularly
varying with index $\delta$ and
%
%e3.4 ###
\begin{equation} \label{hyp_surv2}
\lim_{\xi\to+ \infty}\frac{\psi_r(-\log\bar{F}_r(\mathrm{e}^\xi))}{\xi} \geq 1,
\end{equation}
then
%
%e3.5 ###
\begin{equation} \label{result2_thm_result_gen}
\limsup_{n} \frac{ \log\rho_n}{ \psi_r(\log n)} \leq  \biggl( 1+\frac
{r}{d}  \biggr)^{\delta}.
\end{equation}
If $-\log\bar{F}_r(\mathrm{e}^{\centerdot})$ has regular variation of index $
1/ \delta$ then (\ref{result2_thm_result_gen}) holds with $\psi_r =
(-\log\bar{F}_r(\mathrm{e}^{\centerdot}))^{\leftarrow}$.
\item[(b)] If $\psi_r$ is a non-decreasing function, regularly
varying with index $\delta$ and
%
%e3.6 ###
\begin{equation} \label{hyp_surv1}
\lim_{\xi\to+\infty}\frac{\psi_r(-\log\bar{F}_r(\xi))}{\xi} \geq1,
\end{equation}
then
%
%e3.7 ###
\begin{equation} \label{result1_thm_result_gen}
\limsup_{n } \frac{\rho_n}{ \psi_r(\log n)} \leq c_{r,d}  \biggl( 1+\frac
{r}{d}  \biggr)^{\delta},
\end{equation}
where $c_{r,d} =1$ if $d=1, r \geq1$ and $\mathbf{(G_r)}$ holds and
$c_{r,d}=2$ otherwise.
In particular, if $-\log\bar{F}_r$ has regular variation with index $
1/ \delta$, then (\ref{result1_thm_result_gen}) holds with $\psi_r =
(-\log\bar{F}_r)^{\leftarrow}$.
\end{longlist}
\end{theo}

\textit{Further comments on the choice of $\psi_{r}$.}
As we will show further on, claim (a) is devoted to distributions
with polynomial tails whereas claim (b) will be applied to
distributions with hyper-exponential tails.
Note that for distributions with exponential tails, the function $\psi
_r$ in (b) can be chosen independently of $r$ (see the proof of
Corollary~\ref{cor_princip_limsup}).
Also note that if $-\log\bar{F}_r$ (resp., $-\log\bar
{F}_r(\mathrm{e}^{\centerdot})$) is measurable, locally bounded and regularly
varying with index $1/ \delta, \delta>0$, then its generalized inverse
function $\phi_r$ (resp., $\Phi_r$) is measurable increasing to
${+}\infty$, regularly varying with index $\delta$ and $\phi_r(-\log\bar
{F}_r(x)) = x + \mathrm{o}(x)$ (resp., $\Phi_r(-\log\bar{F}_r(\mathrm{e}^x)) = x
+ \mathrm{o}(x)$). Consequently, inequality (\ref{result1_thm_result_gen})
(resp., (\ref{result2_thm_result_gen})) holds with $\phi_r$
(resp., $\Phi_r$) in place of $\psi_r$. However, $\phi_r$
(resp., $\Phi_r$) is, in general, not easy to compute and the
examples below show that it is often easier to directly exhibit a
function $\psi_r$ satisfying the announced hypotheses without inducing
any asymptotic loss of accuracy.

The above theorem is a consequence of the following more abstract
result, which connects $\rho_n$ and the generalized functions $\bar F_r$.

\begin{theo} \label{prop_princip_limsup}
Let $r>0$ and let $X \in L^{r}(\mathbb{P})$ with a distribution $P$
having an unbounded support and satisfying $\mathbf{(H)}$. Let $(\alpha
_n)_{n \geq1}$ be an $L^r(P)$-optimal sequence of $n$-quantizers. Then,
%
%e3.8 ###
\begin{equation} \label{eqthm2.6}
\lim_{\varepsilon\downarrow0} \liminf_{n}  \biggl( n^{1+  {r}/{d}}
\bar{F}_r \biggl (\frac{\rho_n}{c_{r,d}+\varepsilon}  \biggr)  \biggr) \geq
C_{r,d} \in(0, \infty),
\end{equation}
where $c_{r,d}$ is defined in Theorem~\ref{thm_princip_limsup}.
\end{theo}

We will temporarily admit this result to prove Theorem~\ref{thm_princip_limsup}.

\begin{pf*}{Proof of Theorem~\ref{thm_princip_limsup}}
(a) It follows from (\ref{eqthm2.6}) that, for every $\varepsilon
>0$, there is a positive real constant $C_{r,d,\varepsilon}$
such that $ n^{-\fracb{d+r}{d}} C_{r,d,\varepsilon} \leq\bar{F}_r
(\frac{\rho_n}{c_{r,d} + \varepsilon}  )$. Therefore, one has
\[
\frac{r+d}{d} \log n-\log(C_{r,d,\varepsilon}) \geq-\log\bar{F}_r
\biggl(\frac{\rho_n}{c_{r,d} + \varepsilon}  \biggr).
\]
Combining the fact that $\psi_r$ is non-decreasing with assumption
(\ref{hyp_surv2}) yields
\begin{eqnarray*}
\psi_r  \biggl( \frac{r+d}{d} \log n-\log(C_{r,d,\varepsilon})  \biggr) &
\geq& \psi_r  \biggl(-\log\bar{F}_r  \biggl(\frac{\rho_n}{c_{r,d} +
\varepsilon}  \biggr)  \biggr) \\
& \geq& \log\rho_n - \log(c_{r,d} + \varepsilon)+ \mathrm{o}(\log\rho_n ).
\end{eqnarray*}
Moreover, dividing by $\psi_r(\log n)$ (which is positive for large
enough $n$) yields
\[
\frac{ \log\rho_n}{\psi_r(\log n)} \leq\biggl ( 1-\frac{\log
(c_{r,d}+\varepsilon)}{ \log\rho_n} + \frac{\mathrm{o}(\log\rho_n)}{\log\rho
_n}  \biggr)^{-1} \frac{\psi_r  ( (\fracb{r+d}{d}) \log n-\log
(C_{r,d,\varepsilon})  )}{\psi_r(\log n)}.
\]
Owing to the regularly varying hypothesis on $\psi_r$ and the fact that
$\lim_{n} \rho_n = +\infty$ (which follows from Proposition
\ref{proplimitRM}), we have
\[
\limsup_{n} \frac{ \log\rho_n}{ \psi_r(\log n)} \leq  \biggl(1+ \frac
{r}{d}  \biggr)^{\delta}.
\]

(b) As previously, one derives from (\ref{hyp_surv1}) and from the
non-decreasing hypothesis on~$\psi_r$ that
\begin{eqnarray*}
\psi_r  \biggl( \frac{r+d}{d} \log n-\log(C_{r,d,\varepsilon})  \biggr) &
\geq& \psi_r \biggl (-\log\bar{F}_r \biggl (\frac{\rho_n}{c_{r,d} +
\varepsilon}  \biggr)  \biggr) \\
& \geq& \frac{\rho_n}{c_{r,d} + \varepsilon} + \mathrm{o}(\rho_n ).
\end{eqnarray*}
It follows that
\[
\frac{\rho_n}{\psi_r(\log n)} \leq(c_{r,d}+\varepsilon)  \biggl(1+ \frac
{\mathrm{o}(\rho_n)}{\rho_n}  \biggr)^{-1} \frac{\psi_r  ( (\fracb{r+d}{d}) \log
n-\log(C_{r,d,\varepsilon})  )}{\psi_r(\log n)}.
\]
The regularly varying hypothesis on $\psi_r$ and the fact that
$\lim_{n} \rho_n =+\infty$ yields
\[
\forall\varepsilon>0 \qquad\limsup_{n} \frac{\rho_n}{\psi_r(\log n)}
\leq (c_{r,d}+\varepsilon) \biggl (\frac{r+d}{d}  \biggr)^{\delta}.
\]
The result follows by letting $\varepsilon\rightarrow0.$
\end{pf*}
%
%------------------END PROOF OF THM---------------

Now we pass to the proof of Theorem~\ref{prop_princip_limsup}, which is
based on the following two lemmas.

\begin{lem} \label{lem1_princip_limsup} Let $r>0$ and let $X \in
L^r(\mathbb{P})$ with a distribution $P$ on $\mathbb{R}^d$ having an
unbounded support. Let $(\alpha_n)_{n \geq1}$ be a sequence of
$n$-quantizers, such that $\mathbb{E} d(X,\alpha_n)^r\to0$. Then,
%
%e3.9 ###
\begin{equation} \label{equaconjecture}
\forall \varepsilon>0,  \exists  n_{\varepsilon}  \mbox{ such
that } \forall n\geq n_{\varepsilon}, \forall a \in M(\alpha_n),
\forall y \in C_a(\alpha_n)  \qquad\vert y \vert\geq\frac{\rho
_n}{c_{r,d}+\varepsilon},
\end{equation}
where $c_{r,d}$ is defined in Theorem~\ref{thm_princip_limsup}.
\end{lem}

\begin{pf}
 \textit{Step} 1. Let $r>0$ and let $d \geq
1$. Since $\mathbb{E} d(X,\alpha_n)^r\to0$ as $n \rightarrow+\infty
$, the following asymptotic density property of $(\alpha_n)$ in the
support of $P$ holds:
%
%e3.10 ###
\begin{equation}
\forall\varepsilon>0,  \forall  x \in\operatorname{supp}(P),   \exists
n_{\varepsilon,x} \in\mathbb{N},  \forall n \geq n_{\varepsilon
,x}  \qquad B(x,\varepsilon) \cap\alpha_n \not= \varnothing.
\end{equation}
Otherwise, there exists $x \in\operatorname{supp}(P), \varepsilon>0$ and a
subsequence $(\alpha_{n_k})_{k \geq1}$ so that $ \forall k \geq1$,
$B(x,\varepsilon) \cap\alpha_{n_k} = \varnothing$. Then, for every $k
\geq1$,
\[
\Vert d(X,\alpha_{n_k})\Vert_r \geq\bigl\Vert d(X,\alpha_{n_k})  \mathbf
{1} _{X \in B(x, \varepsilon/2)} \bigr\Vert_r \geq\frac{\varepsilon}{2}
P\bigl(B(x,\varepsilon/2)\bigr)^{1/r} >0,
\]
which contradicts the fact that $\Vert d(X,\alpha_{n})\Vert_r
\rightarrow0$ as $n \rightarrow+\infty$.

Assume first that $0 \in\operatorname{supp}(P)$. Let $ \varepsilon>0$ and $a
\in M(\alpha_n)$. There exists an $N_1  \in\mathbb{N}$ such that
$B(0,\varepsilon) \cap\alpha_n \not= \varnothing$ for every $n\geq N_1$.
Now $\rho_n \rightarrow+\infty$ implies the existence of $N'_1 \in
\mathbb{N}$, $N'_1\ge N_1$ such that $B(0,\varepsilon) \cap(\alpha_n
\backslash M(\alpha_n)) \not= \varnothing$ for $n \geq N'_1$.

Let $n \geq N'_{1}$ and let $b\in B(0,\varepsilon) \cap(\alpha_n
\backslash M(\alpha_n)) $. For every $y \in C_a(\alpha_n)$, we have $
\vert y -b \vert^2 \geq \vert y -a \vert^2$, so that
\[
2(y|a-b) \geq \vert a \vert^2 -\vert b \vert^2 = \rho_n^2 - \vert b
\vert^2\ge0.
\]
Now, if $\vert y \vert\vert a-b \vert\geq (y|a-b$), then,
\[
\vert y \vert\vert a-b \vert\geq\frac{(\rho_n+\vert b \vert)(\rho
_n-\vert b \vert)}{2}.
\]
Moreover, $0<\vert a-b \vert\leq\vert a \vert+ \vert b \vert = \rho
_n+\vert b \vert$. One finally gets
\[
\vert y \vert\geq\frac{\rho_n-\vert b \vert}{2} \geq\frac{\rho
_n-\varepsilon}{2}.
\]
Since $\rho_n \rightarrow+\infty$, then $\vert y \vert\geq\frac{\rho
_n}{2+\varepsilon}$ as soon as $n \geq\max(N'_{1},N_{2})$, with
$N_{2}$ such that $\rho_{N_{2}} \geq2+\varepsilon$.

If $0 \notin\operatorname{supp}(P)$, we show likewise that $\vert y \vert\geq
\frac{\rho_n - \vert x_0 \vert- \varepsilon}{2}$, where $x_0 \in
\operatorname{supp}(P)$ is fixed. This implies the announced result since $\rho
_n \rightarrow+\infty$.

 \textit{Step} 2. Suppose that $d=1, r \geq1$ and ${(\mathbf{G_r})}$
holds. First, we use the well-known fact (see, e.g., \cite{GraLus},
Lemma 4.10 or \cite{Pag}, Proposition 9) that the $L^r$-distortion function
\[
\alpha= (\alpha_1,\ldots ,\alpha_n) \longmapsto D_{n,r}^X (\alpha) =
\mathbb{E} \Bigl( \min_{i=1,\ldots ,n} \vert X - \alpha_i \vert^r  \Bigr)
\]
is differentiable at any codebook $\alpha \in(\mathbb{R}^d)^n$ having
pairwise distinct components and that
%
%e3.11 ###
\begin{equation} \label{station.}
\nabla D_{n,r}^X(\alpha) = r  \biggl( \int_{C_i(\alpha)} (\alpha_i - u)
\vert u - \alpha_i\vert^{r-2} f(u)\,\mathrm{d}u  \biggr)_{1 \leq i \leq n}.
\end{equation}
An optimal $L^r$-quantizer $ \alpha=\{\alpha_1, \ldots ,\alpha_n \}$ at
level $n$ for $P = f \lambda_1$ has full size $n$ so that
%
%e3.12 ###
\begin{equation} \label{EqStationnar}
\nabla D_{n,r}^X(\alpha) =0.
\end{equation}
Note that for any (ordered) quantizer $\alpha_n = \{x_{1}^{(n)},\ldots
,x_{n}^{(n)}\}$, $x_{1}^{(n)}  < \cdots < x_{n}^{(n)}$ at level
$n$, its Voronoi partition is given by
\begin{eqnarray*}
C_{1}(\alpha_n)&=&\bigl(-\infty,x_{ {1}/{2}}^{(n)}\bigr], \qquad  C_n(\alpha
_n)=\bigl(x_{n- {1}/{2}}^{(n)},+\infty\bigr),  \\
  C_i(\alpha_n) &=&\bigl
(x_{i-
{1}/{2}}^{(n)},x_{i+ {1}/{2}}^{(n)}\bigr], \qquad  i=2,\ldots ,n-1,
\end{eqnarray*}
with $ x_{i\pm {1}/{2}}^{(n)} = \frac{x_{i}^{(n)}+x_{i\pm
1}^{(n)}}{2} $. We will focus on the one-sided setting by considering
\[
\rho_n = \rho_n^{+}:= \max\{x, x \in\alpha\}.
\]
All results on $ \rho_n^{-}:= \max\{-x, x \in\alpha\}$ follow by
considering $-X$ instead of $X$. Finally, one will conclude by noting
that the bi-sided radius is given by $\rho_n =\max(\rho^+_n,\rho^-_n)$.

Let $\alpha_n = \{x_1^{(n)},\ldots ,x_n^{(n)} \}$ with $x_1^{(n)}
<  \cdots < x_n^{(n)}$ and suppose that (up to a subsequence) $\frac
{x_{n-1}^{(n)}}{x_{n}^{(n)}} \rightarrow \rho<1$.

Let $\varepsilon>0$ such that $\rho+ \varepsilon<1$. We have for
large enough $n$, $ \frac{x_{n-1}^{(n)}}{x_{n}^{(n)}} < \rho+
\varepsilon<1$
or, equivalently,
%
%e3.13 ###
\begin{equation} \label{equa_ref1}
\frac{x_{n-1}^{(n)}+x_n^{(n)}}{2} < x_n^{(n)} \frac{1+\rho+\varepsilon}{2}.
\end{equation}
Let $\rho'$ be such that $0<\rho'< \frac{1-(\rho+\varepsilon)}{2}$,
that is, $\frac{1+\rho+\varepsilon}{2}< 1-\rho'<1$. It follows from
(\ref{equa_ref1}) that
%
%e3.14 ###
\begin{eqnarray} \label{IneqFromStation}
 \int_{\fracb{x_{n-1}^{(n)} +x_{n}^{(n)}}{2}}^{x_n^{(n)}}  \biggl( 1-\frac
{u}{x_n^{(n)}} \biggr)^{r-1} f(u)\,\mathrm{d}u
& \geq&
 \int_{ {x_n^{(n)}(1+\rho+\varepsilon)}/{2}}^{x_n^{(n)}(1-\rho')}
 \biggl( 1-\frac{u}{x_n^{(n)}} \biggr)^{r-1} f(u)\,\mathrm{d}u \nonumber\\
& \geq& (\rho')^{r-1} \int_{ {x_n^{(n)}(1+\rho+\varepsilon
)}/{2}}^{x_n^{(n)}(1-\rho')} f(u)\,\mathrm{d}u  \\
& \geq& \rho'' x_n^{(n)} f(c_n)\nonumber
\end{eqnarray}
with $\rho'' = (\rho')^{r-1}(\frac{1}{2}-\rho' - \frac{\rho+\varepsilon
}{2}) >0$ and $c_n \in(x_n^{(n)}(1+\rho+\varepsilon
)/2,x_n^{(n)}(1-\rho'))$. On the other hand, since we have
\[
\frac{1}{x_n^{(n)} f(x_n^{(n)})} \int_{x_n^{(n)}}^{+\infty}  \biggl( \frac
{u}{x_n^{(n)}} -1  \biggr)^{r-1}f(u)\,\mathrm{d}u = \int_{1}^{+\infty} (u-1)^{r-1}
\frac{f(u x_n^{(n)})}{f(x_n^{(n)})}\,\mathrm{d}u,
\]
it follows from assumption $\mathbf{(G_r)}$ that
\[
\lim_{n} \frac{1}{x_n^{(n)} f(x_n^{(n)})} \int_{x_n^{(n)}}^{+\infty}
\biggl ( \frac{u}{x_n^{(n)}} -1  \biggr)^{r-1}f(u)\,\mathrm{d}u =0.
\]
Consequently, for large enough $n$,
\[
\frac{1}{x_n^{(n)} f(x_n^{(n)})} \int_{x_n^{(n)}}^{+\infty}  \biggl( \frac
{u}{x_n^{(n)}} -1  \biggr)^{r-1}f(u)\,\mathrm{d}u < \rho''
\]
so that using (\ref{IneqFromStation}) and the fact that $f$ is
non-increasing in $[A,+\infty)$ and $A<c_n < x_n^{(n)}$ for large
enough $n$, one gets
\begin{eqnarray*}
 \int_{x_n^{(n)}}^{+\infty}  \biggl( \frac{u}{x_n^{(n)}} -1
\biggr)^{r-1}f(u)\,\mathrm{d}u &<& \rho'' x_n^{(n)} f\bigl(x_n^{(n)}\bigr)\\
 &\leq& \rho'' x_n^{(n)}
f(c_n) \leq  \int_{\fracb{x_{n-1}^{(n)} +x_{n}^{(n)}}{2}}^{x_n^{(n)}}
\biggl ( 1-\frac{u}{x_n^{(n)}} \biggr)^{r-1} f(u)\,\mathrm{d}u.
\end{eqnarray*}
This leads to a contradiction since the $L^r$-stationary equation (\ref
{EqStationnar}) implies in particular
\[
\int_{\fracb{x_{n-1}^{(n)} +x_{n}^{(n)}}{2}}^{x_n^{(n)}}  \biggl( 1-\frac
{u}{x_n^{(n)}} \biggr)^{r-1} f(u)\,\mathrm{d}u =  \int_{x_n^{(n)}}^{+\infty}  \biggl(
\frac{u}{x_n^{(n)}} -1 \biggr)^{r-1} f(u)\,\mathrm{d}u.
\]
We therefore have shown that $ \lim_{n} \frac
{x_{n}^{(n)}}{x_{n-1}^{(n)}} =1$.
It follows that
\[
\forall\varepsilon>0, \exists n_{\varepsilon}  \mbox{ such that }
\forall n \geq n_{\varepsilon} \qquad  x_{n}^{(n)} < (1+\varepsilon) x_{n-1}^{(n)}.
\]
Thus, one completes the proof by noting that
\[
\forall y \in C_a(\alpha_n), a \in M(\alpha_n) \qquad   \rho_n =
x_{n}^{(n)} < (1+\varepsilon) x_{n-1}^{(n)} < (1+\varepsilon) y.
\]
\upqed
\end{pf}

\begin{lem} \label{lem_minor_diff} Let $r>0$ and let $X \in
L^{r}(\mathbb{P})$ with distribution $P$ satisfying $\mathbf{(H)}$. Let
$(\alpha_n)_{n \geq1}$ be a sequence of $L^r$-optimal $n$-quantizers
of the distribution $P$. Then for large enough $n$,
%
%e3.15 ###
\begin{equation}
e^r_{n,r}(X) - e^r_{n+1,r}(X) \geq C_{r,d} n^{-\fracb{r+d}{d}},
\end{equation}
 {with}
%e3.16 ###
\begin{equation} \label{def_Const}
  C_{r,d} = \frac{r }{2^{(r+d)}(d+r)}  \biggl( \frac
{d}{d+r} \biggr)^{d/r} \frac{\varepsilon_0}{1+\varepsilon_0} Q_{d+r}\bigl(U\bigl(\bar
{B}(x_0, r_0/2)\bigr)\bigr),
\end{equation}
where $U(\bar{B}(x_0,\frac{r_0}{2}))$ stands for the uniform
distribution on the closed ball $\bar{B}(x_0,\frac{r_0}{2})$, the
constants $\varepsilon_{0}$, $x_{0}$, $r_{0}$ come from assumption
$\mathbf{(H)}$ and $Q_{d+r}$ is defined by (\ref{QrP}) in Zador's theorem.
\end{lem}

\begin{pf}
\textit{Step} 1. Let $y \in\mathbb
{R}^d$. We temporarily set $ \delta_n=d(y,\alpha_n)$ and may assume
$\delta_n>0$. Following the lines of the proof of Theorem 2 in \cite
{GraLusPag}, we have for every $x \in B(y,\delta_n/2)$ and $a \in
\alpha_n$,
\[
\vert x -a \vert\geq\vert y-a \vert-\vert x-a \vert\geq\delta_n/2
\]
and hence\vspace*{-1.5pt}
\[
d(x,\alpha_n) \geq\delta_n/2 \geq\vert x-y \vert, \qquad x \in
B(y,\delta_n/2).
\]
It follows, by setting $\beta_n = \alpha_n \cup\{ y\}$, that
$d(x,\alpha_n)\ge d(x,\beta_n)$ and $ d(x,\beta_n) = \vert x -y \vert,
  x \in B(y,\delta_n/2).$ Consequently for every $b \in(0,1/2)$,
\begin{eqnarray*}
e^r_{n,r}(X) -e^r_{n+1,r}(X)
& \geq& \int_{B(y,\delta_n b)}\bigl(d(x,\alpha_n)^r-d(x,\beta_n)^r\bigr) P(\mathrm{d}x)
\\
& = &  \int_{B(y,\delta_n b)} \bigl(d(x,\alpha_n)^r- \vert x-y \vert^r\bigr)
P(\mathrm{d}x) \\
& \geq& \int_{B(y,\delta_n b)} \bigl((\delta_n/2)^r -( \delta_n
b)^r\bigr)P(\mathrm{d}x) \\
& = & (2^{-r} -b^r) \delta_n^r P(B(y,\delta_n b)).
\end{eqnarray*}

\textit{Step} 2. This step is the core of our proof. Let $x_0$
and $r_0$ be as in $\mathbf{(H)}$. For every $y \in\bar{B}(x_0,\frac
{r_0}{2})$,
\begin{eqnarray*}
e^r_{n,r}(X) - e^r_{n+1,r}(X) & \geq& (2^{-r} - b^{r}) \delta_n^r
P\biggl(B\biggl(y, \min\biggl(b \delta_n, \frac{r_0}{2} \biggr)\biggr)\biggr) \\
& \geq& (2^{-r} - b^{r}) \delta_n^r \varepsilon_0 \min \biggl((b \delta
_n)^d, \biggl(\frac {r_0 }2\biggr)^d  \biggr).
\end{eqnarray*}
We know from \cite{DelGraLusPag} that, as soon as $d(x,\alpha_n)\to0$
as $n\to\infty$ in $L^r(P)$, the convergence will hold uniformly on
compact sets as well. In particular, we have
\[
\sup_{y\in\bar{B}(x_0, {r_0}/{2})} d(y,\alpha_n) \rightarrow0
\]
so that there exists $N(x_0,r_0) \in\mathbb N$ such that for every
$n \ge N(x_0,r_0)$,
\[
\sup_{y\in\bar{B}(x_0, {r_0}/{2})} d(y,\alpha_n) \leq\frac{r_0}{2}.
\]
Consequently\vspace*{-1.5pt}
\[
e^r_{n,r}(X) - e^r_{n+1,r}(X) \geq(2^{-r} - b^{r}) b^d d(y,\alpha
_n)^{d+r} \varepsilon_0 \mathbf{1} _{\{y\in\bar{B}(x_0,
{r_0}/{2})\}}.
\]
It follows that
\begin{eqnarray*}
e^r_{n,r}(X) - e^r_{n+1,r}(X)
& \geq& (2^{-r} - b^{r}) \varepsilon_0 b^d \int_{\bar{B}(x_0,
{r_0}/{2})} d(y,\alpha_n)^{d+r}\frac{\lambda_d(\mathrm{d}y)}{\lambda_d(\bar
{B}(x_0, {r_0}/{2}))} \\
& \geq& (2^{-r} - b^{r})b^d  \varepsilon_0 \lambda_d\bigl(\bar
{B}(x_0,r_0/2)\bigr) e_{n,r+d}^{r+d}\bigl(U\bigl(\bar{B}(x_0,r_0/2)\bigr)\bigr),
\end{eqnarray*}
where we used in the last inequality the fact that $\alpha_n$ is
suboptimal for the uniform distribution over $\bar{B}(x_0,\frac
{r_0}{2}).$ As a consequence,
\[
e^r_{n,r}(X) - e^r_{n+1,r}(X) \geq (2^{-r} - b^{r})   b^d \varepsilon
_0  e_{n,r+d}^{r+d}\bigl(U\bigl(\bar{B}(x_0,r_0/2)\bigr)\bigr).
\]
Finally, one completes the proof by noting that, for large enough $n
\geq N(x_0,r_0)$,
\[
e^r_{n,r}(X) - e^r_{n+1,r}(X) \geq\sup_{b\in(0,1/2)} \bigl((2^{-r} -
b^{r})b^d\bigr) \frac{\varepsilon_0}{1+\varepsilon_0}  Q_{d+r}\bigl(U\bigl(\bar
{B}(x_0,r_0/2)\bigr)\bigr) n^{-\fracb{d+r}{d}}.
\]
\upqed
\end{pf}

Now we are in position to complete the proof of Theorem~\ref
{prop_princip_limsup}.
%-------------------PROOF OF THE THM----------------------
%
\begin{pf*}{Proof of Theorem \ref{prop_princip_limsup}}
Let $a \in M(\alpha_n)$ and $\varepsilon>0$. We have,
\[
e^r_{n-1,r}(X) =\mathbb{E} \vert X - \widehat{X}^{\alpha_{n-1 }} \vert
^r \leq\mathbb{E} \bigl\vert X - \widehat{X}^{\alpha_{n} \backslash\{a\} }
\bigr\vert^r
\]
since $\alpha_{n-1}$ is $L^r$-optimal at level $n-1$. Hence
\begin{eqnarray*}
\mathbb{E} \bigl\vert X - \widehat{X}^{\alpha_{n} \backslash\{ a \} } \bigr\vert
^r & = &
\mathbb{E}  \bigl( \vert X - \widehat{X}^{\alpha_{n} } \vert^r  \mathbf
{1} _{ \{ X \in C_a^{c}(\alpha_{n}) \}}  \bigr) +
\mathbb{E}  \Bigl( \min_{b \in\alpha_{n} \backslash\{a \}} \vert X - b
\vert^r  \mathbf{1} _{ \{ X \in C_a(\alpha_{n}) \}}  \Bigr) \\
& \leq & e^r_{n,r}(X) + \mathbb{E}  \Bigl( \min_{b \in\alpha_{n}
\backslash\{ a \} } ( \vert X \vert+ \vert b \vert)^r \mathbf
{1}_{ \{ X \in C_a(\alpha_{n}) \}}  \Bigr).
\end{eqnarray*}
It follows from Lemma \ref{lem1_princip_limsup} that, for every
$\varepsilon>0$, there exists $n_{\varepsilon} \in\mathbb{N}$ such
that for every $n \geq n_{\varepsilon}$, $\vert X \vert> \frac{\rho
_n}{c_{r,d}+\varepsilon}$, on the event $\{ X \in C_a(\alpha_{n}) \}$.
Consequently, for all $ b \in\alpha_{n} \backslash\{ a \}, \vert b
\vert\leq\vert a \vert= \rho_n < (c_{r,d}+\varepsilon) \vert X \vert
$. Hence,
\[
e^r_{n-1,r}(X) - e^r_{n,r}(X) \leq(c_{r,d}+1+\varepsilon)^r \mathbb
{E}  \bigl(\vert X \vert^r \mathbf{1}_{\{ \vert X \vert>  {\rho
_n}/({c_{r,d}+\varepsilon}) \}}  \bigr).
\]
Lemma \ref{lem_minor_diff} yields for large enough $n$ (since
$(n-1)^{-\fracb{r+d}{d}} \sim n^{-\fracb{r+d}{d}}$ as $n \rightarrow
+\infty$),
\[
(1+\varepsilon)^{-1} C_{r,d} n^{-\fracb{r+d}{d}} \leq
(c_{r,d}+1+\varepsilon)^r \mathbb{E}  \bigl(\vert X \vert^r \mathbf{1}_{\{ \vert X \vert>  {\rho_n}/({c_{r,d}+\varepsilon}) \}}  \bigr)
\]
so that for every $\varepsilon>0$,
\[
\liminf_{n }  \biggl( n^{\fracb{r+d}{d}} \bar{F}_r  \biggl(\frac{\rho
_n}{c_{r,d}+\varepsilon}  \biggr)  \biggr) \geq\frac{C_{r,d}}{
(c_{r,d}+1+\varepsilon)^r (1+\varepsilon) }.
\]
Letting $\varepsilon\rightarrow0$ yields the statement (\ref{eqthm2.6}).
\end{pf*}

%s3.2 ###
\subsection{Applications to distributions with polynomial and
hyper-exponential tails} \label{Upperexplicit}
We next give an explicit
asymptotic upper bound for the convergence rate of the maximal radius
sequence by making the function $\psi_r$ explicit. These bounds are
derived in terms of the rate of decay of the generalized survival
function $\bar{F}_r$.

\begin{prop} \label{cor_princip_limsup} Let $r>0$ and let $X \in
L^{r+}(\mathbb{P})$ with distribution $P$ having an unbounded support
and satisfying $(\mathbf{H})$. Suppose that $(\alpha_n)_{n \geq1}$ is
an $L^r$-optimal sequence of $n$-quantizers for~$X$.
\begin{longlist}[(b)]
\item[(a)] \textup{Polynomial tail.} Set
%
%e3.17 ###
\begin{equation} \label{Assum_Cor_pol_tail1}
\zeta^{\star}=\sup \Bigl\{ \zeta>0, \limsup_{\xi\rightarrow+\infty} \xi
^{\zeta-r} \bar{F}_{r}(\xi) <+\infty  \Bigr\} = \sup \{ \zeta>r,
\mathbb{E} \vert X \vert^{\zeta} <+\infty  \}.
\end{equation}
Then $\zeta^{\star} \in(r,+\infty]$ and
%
%e3.18 ###
\begin{equation} \label{Ineq_Cor_pol_tail1}
\limsup_{n} \frac{\log\rho_n}{\log n} \leq\frac{1}{\zeta^{\star}-r}
\frac{r+d}{d}.
\end{equation}
\item[(b)] \textup{Hyper-exponential tail.} Assume there exists $\kappa>0$
such that $\mathrm{e}^{\vert X \vert^{\kappa}} \in L^{0+}(\mathbb{P})$. Set
%
%e3.19 ###
\begin{equation} \label{Assum_Cor_Exp_tail1}
{\theta}^{\star} = \sup \Bigl\{ \theta>0, \limsup_{\xi\rightarrow
+\infty} \mathrm{e}^{\theta \xi^{\kappa}} \bar{F}_r(\xi) < +\infty \Bigr\} =
\sup \bigl\{ \theta>0, \mathbb{E} \mathrm{e}^{\theta\vert X \vert^{\kappa}} <
+\infty \bigr\}.
\end{equation}
Then $\theta^{\star} \in(0,+\infty]$ and
%
%e3.20 ###
\begin{equation} \label{Ineq_Cor_Exp_tail1}
\limsup_{n} \frac{\rho_n}{ (\log n  )^{1/ \kappa}} \leq c_{r,d}
\biggl ( \frac{r+d}{ d \theta^{\star} }  \biggr)^{1/ \kappa}.
\end{equation}
\end{longlist}
\end{prop}

\begin{rem}   If $X \in\bigcap_{r>0} L^r(\mathbb{P}),$
then $\zeta^{\star}=+\infty$ and, consequently, $ \lim_{n\rightarrow+\infty} \frac{\log\rho_n}{\log n} =0$. This
confirms that this asymptotics is not the significant one for
distributions with hyper-exponential tails.\vspace*{-3pt}
\end{rem}

\begin{pf*}{Proof of Proposition~\ref{cor_princip_limsup}}
The equalities in (\ref{Assum_Cor_Exp_tail1}) and (\ref
{Assum_Cor_pol_tail1}) are elementary.
\begin{longlist}[(b)]
\item[(a)] Let $\zeta \in(r,\zeta^{\star})$. We have
\begin{eqnarray*}
\mathbb{E}  \bigl(\vert X \vert^r \mathbf{1}_{\{ \vert X \vert> \xi\}
}  \bigr)
& = & \mathbb{E}  \bigl(\vert X \vert^r \mathbf{1}_{\{ 1< \xi^{-\zeta
+r} \vert X \vert^{\zeta-r} \}}  \bigr) \\[-2pt]
& \leq& \xi^{-\zeta+r} \mathbb{E} \vert X \vert^{\zeta}.
\end{eqnarray*}
Then $ -\log\bar{F}_r(\xi) \geq(\zeta-r) \log\xi+ C$, $C \in
\mathbb{R}$, so that by setting $\psi_r(\xi) = \frac{\xi}{\zeta- r}$,
it follows from Theorem~\ref{thm_princip_limsup}(a) that
\[
\limsup_{n} \frac{\log\rho_n}{\log n} \leq\frac{1}{\zeta-r} \frac{r+d}{d}.
\]
Letting $\zeta$ go to $\zeta^{\star}$ yields the assertion (\ref
{Ineq_Cor_pol_tail1}).
\item[(b)] Let $\theta \in(0,\theta^{\star})$. We have
\[ \label{Ineq_tail}
\mathbb{E}  \bigl(\vert X \vert^r \mathbf{1}_{\{ \vert X \vert> \xi
\}}  \bigr) = \mathbb{E} \bigl (\vert X \vert^r \mathbf{1}_{ \{
\mathrm{e}^{\theta\vert X \vert^{\kappa}} > \mathrm{e}^{\theta\xi^{\kappa}}  \} }
 \bigr) \leq \mathrm{e}^{-\theta\xi^{\kappa}} \mathbb{E} \bigl (\vert X \vert^r
\mathrm{e}^{\theta\vert X \vert^{\kappa}} \bigr).
\]
Now, the right-hand side of this last inequality is finite because if
$\theta' \in(\theta,\theta^{\star})$, there exists a positive
constant $C_{\theta, \theta'} $ such that, for every $\xi \in\mathbb
{R}^d$, $\vert\xi\vert^r \mathrm{e}^{\theta\vert \xi\vert^{\kappa}} \leq1
+ C_{\theta, \theta'} \mathrm{e}^{\theta' \vert \xi\vert^{\kappa}}$. As a~consequence,
\[
-\log\bar{F}_r(\xi) \geq\theta\xi^{\kappa} + C_{\theta,X}, \qquad
C_{\theta,X} \in\mathbb{R}.
\]
Let $\psi_{\theta} (y)=  ( \frac{y}{\theta}  )^{1/ \kappa}$.
As a function of $y$, $\psi_{\theta}$ is continuous increasing to
${+}\infty$, regularly varying with index $\delta=\frac{1}{\kappa}$ and
we have
\[
\psi_{\theta}(-\log\bar{F}_r(\xi)) \geq \biggl(\xi^{\kappa} + \frac
{C_X}{\theta} \biggr)^{1/ \kappa} = \xi+ \mathrm{o}(\xi)\qquad\mbox{as } \xi
\rightarrow+\infty.
\]
It follows from Theorem~\ref{thm_princip_limsup}(b)  that, for every
$\theta \in(0,\theta^{\star})$,
\[
\limsup_{n} \frac{\rho_n}{ (\log n  )^{1/ \kappa}} \leq c_{r,d}
 \biggl(\frac{d+r}{d \theta}  \biggr)^{1/ \kappa}.
\]
Letting $\theta\rightarrow\theta^{\star}$ completes the proof.
\end{longlist}
\upqed\vspace*{-3pt}
\end{pf*}

We now give more explicit results for two wide classes of density
functions in $\mathbb{R}^d$: the distributions with polynomial tails
and hyper-exponential tails which, among others, include the Pareto,
Gaussian, Weibull, gamma and double-sided gamma distributions, respectively.\vspace*{-3pt}

\begin{cor} \label{cor_gen_density}
\textup{(a)} If the density $f$ of $X$ satisfies
%
%e3.21 ###
\begin{equation} \label{gen_density_pol_limsup}
\limsup_{|x|\to+ \infty}  |x|^c f(x)<+\infty\qquad\mbox{for some }
c>r+d,\vadjust{\goodbreak}
\end{equation}
then $X \in L^{r+}(\mathbb{P})$ and
%
%e3.22 ###
\begin{equation} \label{Ineq_Cor_pol_tail_density2}
\zeta^{\star} \ge c-d\quad\mbox{and}\quad \limsup_{n}   \frac
{\log\rho_n}{\log n} \leq\frac{1}{c-d-r} \frac{r+d}{d}.
\end{equation}

\textup{(b)} If the density of $X$ satisfies
%
%e3.23 ###
\begin{equation} \label{gen_density_exp_limsup}
\limsup_{|x|\to+ \infty} \frac{\log f(x)}{|x|^{\kappa}} =-\vartheta
<0\qquad\mbox{for some } \kappa>0,
\end{equation}
then $X \in L^{r+}(\mathbb{P})$ and
%
%e3.24 ###
\begin{equation} \label{Ineq_Cor_Exp_tail_density1}
\theta^{\star}\ge \vartheta \quad\mbox{and}\quad \limsup_{n}
\frac{\rho_n}{ (\log n  )^{1/ \kappa}} \leq\frac
{c_{r,d}}{\vartheta^{1/ \kappa}}  \biggl(1+ \frac{r}{d}  \biggr)^{1/ \kappa}.
\end{equation}
\end{cor}

\begin{pf}
(a)
Let $A$, $B>0$ such that for every $x$ with $|x|\ge B$, $
f(x)\le\frac{A}{|x|^c}$. Then, as soon as $\xi\ge B$,
\begin{eqnarray*}
\bar F_r(\xi) = \mathbb{E} \bigl( |X|^r\mathbf{1}_{\{|X|\ge\xi\}}
 \bigr)\le A \int_{\{|x|\ge\xi\}} |x|^r \frac{\mathrm{d}x}{|x|^c}= A  d V_d
 \operatorname{det}(S) \frac{\xi^{ r+d-c}}{r+d-c}, &&
\end{eqnarray*}
where $V_d$ denotes the hyper-volume of the unit Euclidean ball of
$\mathbb{R}^d$ and $|x|^2=  ^t xS x$. As a consequence, for any $\zeta
<c-d$ and any $\xi\ge B$,
\[
\xi^{\zeta-r}\overline F_r(\xi) \le A  d V_d \operatorname{det}(S) \frac{\xi
^{ r+d-c}}{r+d-c}
\]
so that $\overline{\lim}_{\xi\to\infty}
\xi^{ \zeta-r} \bar
F_r(\xi) =0, $ that is, $\zeta^{\star}\ge c-d$ by Proposition~\ref
{cor_princip_limsup}(a).

 (b) It follows from the assumption that, for every $\eta\in
(0,\vartheta/3)$, there exists $B>0$ such that, for every $x$ with
$|x|\ge B$, $ f(x)\le \mathrm{e}^{-(\vartheta-\eta)|x|^{\kappa}}$.
Hence, as soon as $\xi\ge B$,
\begin{eqnarray*}
\bar F_r(\xi) &=& \mathbb{E} \bigl( |X|^r\mathbf{1}_{\{|X|\ge\xi\}
} \bigr)\\
&\le&  \int_{\{|x|\ge\xi \}} |x|^r \mathrm{e}^{-(\vartheta-\eta)|x|^{\kappa
}}\,\mathrm{d}x\\
&= & d  V_d \operatorname{det}(S) \int_{\{ u\ge\xi\}}u^{r+d-1}
\mathrm{e}^{-(\vartheta-\eta) u^{\kappa}}\,\mathrm{d}u
\end{eqnarray*}
so that
\[
  \mathrm{e}^{(\vartheta-3\eta)\xi^{\kappa}}\bar F_r(\xi
)  \le d  V_d  \operatorname{det}(S)  \mathrm{e}^{-\eta\xi^{\kappa}}\int_{\{u\ge B \}}
u^{r+d-1}\mathrm{e}^{-\eta u^{\kappa}}\,\mathrm{d}u.
\]
Consequently, $\theta^{\star}\ge\vartheta-3\eta$ and letting $\eta$ go
to $0$ shows that $\theta^{\star}\ge\vartheta$, which completes the proof.
\noqed\mbox{\hfill}\qed
\end{pf}

%--------------------------------- LOWER ESTIMATE
%-------------------------------
%s4 ###
\section{Lower estimate and asymptotic rates}\label{sec4}
In this section we study the asymptotic lower estimate of the maximal
radius sequence $(\rho_n)_{n \geq1}$ induced by an $L^{r}$-optimal
sequence of $n$-quantizers. First we introduce the family of the
$(r,s)$-distributions, which will play a crucial role to obtain the
sharp lower estimate of the maximal radius sequence.

Let $r>0$ and let $ s>r$. Since the $L^r$-norm is increasing, it is
clear that, for every $s \leq r$, any $L^r$-optimal sequence of
quantizers $(\alpha_n)_{n \geq1}$ is \textit{$L^s$-rate optimal}, that
is, $\limsup_{n} n^{1/d} \Vert X - \widehat{X}^{\alpha_n} \Vert_s <
+\infty.$

But if $s>r$ (and $X \in L^s(\mathbb{P})$), this asymptotic rate
optimality usually fails. This is always the case when $s>r+d$ and $X$
has a probability distribution $f$ satisfying $\lambda_d(f>0) = +\infty
$, as pointed out in \cite{GraLusPag}, Corollaries 3 and 4. It is
established in \cite{Sag} that some linear transformation of the
$L^r$-optimal quantizers $(\alpha_n)$ makes it possible to overcome the
critical exponent $r+d$; that is, one can always construct an
$L^s$-rate-optimal sequence of quantizers up to an affine
transformation of the $L^r$-optimal sequence of quantizers $(\alpha_n)$.

However, there are many (usual) distributions for which $L^s$-rate
optimality does hold for every $s \in[ r,r+d)$. This leads to the
following definition.

\begin{defin} \label{def_(r,s)_dist}
Let $r>0$ and $\nu \in(0,d)$. A random vector $X \in L^{r+}(\mathbb
{P})$ has an \textit{$(r,r+\nu)$-distribution} if any $L^r$-optimal
sequence $(\alpha_n)_{n \geq1}$ is $L^{r+\nu}$-rate optimal, that is,
\[
\limsup_{n} n^{1/d} \Vert X - \widehat{X}^{\alpha_n} \Vert_{r+\nu} <
+\infty.
\]
\end{defin}

Note that if $X$ has an $(r,r+\nu)$-distribution, then $X \in L^{r+\nu
}(\mathbb{P})$. A necessary condition for a distribution $P$ with
density $f$ to have an $(r,r+\nu)$-distribution is (see \cite{GraLusPag}):
%
%e4.1 ###
\begin{equation} \label{necess_rate_opt}
\int_{\mathbb{R}^d} f(x)^{-\fracc{(r+\nu)}{d+r}} P(\mathrm{d}x) < +\infty.
\end{equation}

For $\nu \in(0,d)$, criterions that imply that $X$ has an $(r,r+\nu
)$-distribution have been provided in \cite{GraLusPag}. We mention two
of them below.

\begin{prop}[(Radial tail)] \label{crit_radial}  Let $r>0$ and let $X
\in L^{r+}(\mathbb{P})$ with distribution $P= f \lambda_d$ having an
unbounded support that is the intersection of finitely many half-spaces.
\begin{longlist}[(b)]
\item[(a)] Suppose $f$ has a \textit{radial tail}, that is, there
exists a norm $N(\cdot)$ on $\mathbb{R}^{d}$ and $R_0 \in\mathbb{R}_{+}$
such that
%
%e4.2 ###
\begin{equation}\label{Essradial}
%\begin{tabular}{@{\hspace*{-25pt}}p{338pt}@{}}
\hspace*{-15pt}f = h(N( \cdot))    \mbox{ on } B_{N(\cdot)}(0,R_0)^{c}, \mbox{ where }
h\dvtx
[R_0,+\infty) \rightarrow \mathbb{R}_{+} \mbox{ is a decreasing
function.}
%\end{tabular}
\end{equation}
Let $\nu \in(0,d)$. If
%
%e4.3 ###
\begin{equation} \label{asscor1}
 \int_{\mathbb{R}^d} f(\rho x)^{- ({r+\nu})/({r+d})} P(\mathrm{d}x) < +\infty
\qquad\mbox{for some $\rho>1,$}
\end{equation}
then $X$ has an $(r,r+\nu)$-distribution.

\item[(b)] Assume $d=1$. If $ \operatorname{supp}(P) \subset[A_0, +\infty
)$ for some $A \in\mathbb{R}$, $f_{|(R_0, +\infty)}$ is decreasing
for $R_0 \geq A_0$ and, if, furthermore, assumption
(\ref{asscor1}) holds for some $\rho>1$, then $X$ has an $(r,r+\nu
)$-distribution.
\end{longlist}
\end{prop}

The following proposition works for distributions with non-radial tails.

\begin{prop} \label{Critere2} Let $r>0$ and let $X \in L^{r+}(\mathbb
{P})$ with distribution $P= f \lambda_d$ having a~convex (unbounded)
support. Assume that $f$ satisfies the following \textit{local decay
control} assumption: There exist real numbers $\varepsilon\geq0, \eta
 \in(0,1)$, $M$, $ K >0$ such that
%
%e4.4 ###
\begin{equation}\label{LocalControl}
\forall x,y \in\operatorname{supp}(P),  \vert x \vert\geq M,   \vert y-x
\vert\leq \eta\vert x \vert  \quad \Longrightarrow \quad  f(y) \geq K
f(x)^{1+\varepsilon}.
\end{equation}
Let $\nu \in(0,d)$. If
%
%e4.5 ###
\begin{equation} \label{Hypoth2}
\int_{\mathbb{R}^d} f(x)^{-\fracc{(r+\nu)(1+\varepsilon)}{r+d}} P(\mathrm{d}x) <
+\infty,
\end{equation}
then $X$ has an $(r,r+\nu)$-distribution.
\end{prop}

It follows from Proposition~\ref{crit_radial} that the Gaussian,
Weibull and gamma distributions are all $(r,r+\nu)$-distributions for
every $\nu \in(0,d)$. The Pareto distribution with index $\gamma>r$
has an $(r,r+\nu)$-distribution if and only if $ \nu \in(0,\frac
{\gamma-r}{\gamma+1})$.

More generally, if a distribution $P= f \lambda_d$ is supported by a
convex subset $C$ of $\mathbb{R}^d$ such that
\[
f(x)= \mathrm{e}^{-g(x)^{\kappa}}, \qquad g\dvtx C\to\mathbb{R}_+,  \mbox{ Lipschitz
continuous, }\kappa>0,
\]
or
\[
f(x) = \frac{1}{g(x)^c}, \qquad g\dvtx C\to \mathbb{R}_+,  \mbox{ Lipschitz
continuous, }g\ge\varepsilon_0 \mbox{ on }   B(0,M)^{c},   c>d,
\]
then $P$ satisfies the local decay control criterion (\ref
{LocalControl}) of Proposition~\ref{Critere2} for arbitrarily small
positive $\varepsilon$ and $\varepsilon=0,$ respectively.

Now, suppose that $X$ has an $(r,r+\nu)$-distribution for some $\nu
\in(0,d)$ and set
\[
\nu^{\star}_{X}:= \sup\{\nu>0 \mbox{ s.t. $X$ has an $(r,r+\nu
)$-distribution}\}  \in[0,d].
\]
Note that $X \in L^{r+\nu}(\mathbb{P})$ for every $\nu  \in(0,\nu
^{\star}_{X})$ and that
\[
\{\nu>0 \mbox{ s.t. $X$ has an $(r,r+\nu)$-distribution}\} = (0,\nu
^{\star}_{X}) \mbox{ or } (0,\nu^{\star}_{X}].
\]
When $ \{\nu>0 \mbox{ s.t. $X$ has an $(r,r+\nu)$-distribution}\} =
\varnothing$, we set
\[
\nu^{\star}_{X}=0^+  \qquad \mbox{with the convention }  [0,0^+)=
\{0\}.\vspace*{-2pt}
\]
This convention is consistent with the Zador theorem satisfied by $X
\in L^{r+}(\mathbb{P})$. Note that~$\nu^{\star}_{X}$ may be lower than
$d$, as is the case for the Pareto distribution.

We present below two different approaches to derive the asymptotic
lower bound. The first one is based on tail estimates and involves the
generalized survival functions $\bar{F}_r$ like for the upper estimate.
The second one is based on a new connection with mean random quantization.

%------------------------------ DISTRIBUTION TAIL APPROACH
%----------------------------

%s4.1 ###
\subsection{Distribution tail approach}\vspace*{-2pt}\label{sec4.1}
%s4.1.1 ###
\subsubsection{General results on asymptotic lower bounds}\label{sec4.1.1}
The main result of this section is the theorem below, which connects
the asymptotic lower estimate for $\rho_n$ with the regularly varying
property of ``the'' asymptotic inverse of $-\log\bar{F}_r$ (or one of
its lower bound).

\begin{theo} \label{thm_princip_liminf} Let $r>0$ and let $X \in
L^{r+}(\mathbb{P})$ be an $\mathbb{R}^d$-valued random variable with
distribution $P$ having an unbounded support.
Let $(\alpha_n)_{n \geq1}$ be an $L^r(P)$-optimal sequence of $n$-quantizers.
\begin{longlist}[(b)]
\item[(a)] Let $\nu \in[0,\nu^{\star}_{X})$. If there is a
non-decreasing function $\psi_{r,\nu}$ going to ${+}\infty$ as $x \to
+\infty$, regularly varying with index $\delta$ and satisfying
%
%e4.6 ###
\begin{equation} \label{Hyp_liminf2}
\limsup_{\xi\to+ \infty}\frac{\psi_{r,\nu}(-\log\bar{F}_{r+\nu}(\mathrm{e}^{\xi
}))}{\xi} \leq 1,
\end{equation}
then
%
%e4.7 ###
\begin{equation} \label{eq_thm_gen_liminf2}
\liminf_{n} \frac{\log\rho_n}{\psi_{r,\nu}(\log n)} \geq  \biggl( \frac
{r + \nu}{d}  \biggr)^{\delta}.
\end{equation}
In particular, if $-\log\bar{F}_{r+\nu}(\mathrm{e}^{x})$ has regular variation
with index $1/ \delta$, then (\ref{eq_thm_gen_liminf2}) holds with
$\psi_{r,\nu}(x)= (-\log\bar{F}_{r+\nu}(\mathrm{e}^{x}))^{\leftarrow}$.
\item[(b)] If $\psi$ is a non-decreasing function going to ${+}\infty
$ as $x \rightarrow+\infty$, regularly varying with index $\delta$ and
satisfying
%
%e4.8 ###
\begin{equation} \label{Hyp_liminf1}
\limsup_{\xi\to+ \infty}\frac{ \psi(-\log\bar{F}(\xi))}{\xi} \leq1,
\end{equation}
then
%
%e4.9 ###
\begin{equation} \label{eq_thm_gen_liminf1}
\liminf_{n} \frac{\rho_n}{ \psi(\log n) } \geq  \biggl( \frac{r +\nu
^{\star}_{X}}{d}  \biggr)^{\delta}.
\end{equation}
If $-\log\bar{F}$ has regular variation of index $1/ \delta$, then
(\ref{eq_thm_gen_liminf1}) holds with $\psi= (-\log\bar
{F})^{\leftarrow}$.
\end{longlist}
\end{theo}

Similar to the upper limit, one may note that for distribution with
exponential tails, the function $\psi$ does not depend on $r$ and $\nu$
even if in assumption (\ref{Hyp_liminf1}) we take the generalized
survival function $\bar{F}_{r+\nu}$ instead of the regular survival
function $\bar{F}$. However, for distributions with polynomial tails
like the Pareto distribution, the function $\psi_{r,\nu}$ in~(\ref
{Hyp_liminf2}) may depend on $r$ and consideration of the standard
survival function $\bar{F}$ in place of $\bar{F}_{r+\nu}$ would lead to
a less accurate lower bound.\vadjust{\goodbreak}

As for the upper limit, this result essentially relies on a more
abstract result that connects $\rho_n$ and the (generalized) survival
functions $\bar F_r$.

\begin{theo} \label{prop_gen_liminf}
Let $r>0$ and let $X \in L^{r+}(\mathbb{P})$ be an $\mathbb
{R}^d$-valued random variable with distribution $P$. Let $(\alpha_n)_{n
\geq1}$ be an $L^r(P)$-optimal sequence of $n$-quantizers. For every
$\nu \in[0,\nu^{\star}_{X})$, the following statements hold:
%
%e4.11 ###
%e4.10 ###
\begin{eqnarray}\label{eqthm1.6}
 &&\hspace*{-112pt}\textup{(a)} \quad   \displaystyle \limsup_{n} \sup_{c>0}  \bigl( c^{r+\nu} n^{\fracb{r+\nu
}{d}} \bar{F} ( \rho_n + c )  \bigr) < +\infty.\\
 \label{eqthm1.6.1}
&&\hspace*{-112.5pt} \textup{(b)} \quad   \displaystyle \limsup_{n} \sup_{u >1}  \bigl(  ( 1- 1/u  )^{r+\nu
} n^{\fracb{r+\nu}{d}} \bar{F}_{r+\nu}  ( u \rho_n  )  \bigr) <
+\infty.
\end{eqnarray}
\end{theo}

We temporarily admit this theorem to prove Theorem~\ref{thm_princip_liminf}.

%--------------------PROOF OF THEOREM------------------
%
\begin{pf*}{Proof of Theorem~\ref{thm_princip_liminf}} Let
us focus on (b) (claim (a) is proved in a similar manner by
considering $\bar{F}_{r+\nu}$ instead of $\bar{F}$, for $\nu \in[0,\nu
^{\star}_{X})$). Let
$ \nu \in[0,\nu^{\star}_{X})$. It follows from (\ref{eqthm1.6})
that for large enough $n$,
\[
- \log\bar{F}( \rho_n +c) \geq- \log(C_{\nu,c}) + \frac{r+\nu}{d}
\log n,
\]
where $C_{\nu,c}$ is a positive real constant depending on the indexing
parameters. We derive from the fact that $\psi$ is non-decreasing and
goes to ${+}\infty$ and from assumption (\ref{Hyp_liminf1}) that
\[
\frac{\rho_n}{\psi(\log n)} \geq \biggl( 1+ \frac{c}{\rho_n} + \frac
{\mathrm{o}(\rho_n)}{ \rho_n}  \biggr)^{-1} \frac{\psi( \fracb{r+\nu}{d} \log n -
\log(C_{\nu,c})) }{\psi(\log n)}.
\]
Since $\psi$ is regularly varying with index $\delta$ we have
\[
\forall\nu \in[0,\nu^{\star}_{X}) \qquad \liminf_{n} \frac{\rho
_n}{\psi(\log n)} \geq  \biggl( \frac{r+\nu}{d}  \biggr)^{\delta}.
\]
When $\nu^{\star}_{_X}>0$, letting $\nu\rightarrow\nu^{\star}_{X}$
yields the announced result.
\end{pf*}

\begin{pf*}{Proof of Theorem~\ref{prop_gen_liminf}}
(a)
Let $n\ge1$, let $c>0$ and let $\nu \in[0,\nu^{\star}_{X})$. Then
\[ \label{eq.bad_estim}
\mathbb{E} \vert X - \widehat{X}^{\alpha_n} \vert^{r+\nu} \geq
\mathbb{E}  \Bigl( \min_{a \in\alpha_n}\vert X - a \vert^{r+\nu}  \mathbf{1} _{ \{ \vert X \vert> \rho_n + c \}}  \Bigr).
\]
In the event $\{\vert X \vert> \rho_n + c \}$, we have $\vert X \vert
> \rho_n +c > \rho_n \geq \vert a \vert$ for every $a \in\alpha_n$. Then
%
%e4.12 ###
\begin{eqnarray}\label
{estim_error}\label{result_inter} %\label{result_inter}
n^{\fracb{r+\nu}{d}}\mathbb{E} \vert X - \widehat{X}^{\alpha_n} \vert
^{r+\nu} & \geq&n^{\fracb{r+\nu}{d}}
\mathbb{E}  \Bigl( \min_{a \in\alpha_n}\vert X - a \vert^{r+\nu}
{\mathbf{1}}_{ \{ \vert X \vert> \rho_n + c \}}  \Bigr)  \nonumber \\
& \geq& n^{\fracb{r+\nu}{d}}\mathbb{E}  \Bigl( \min_{a \in\alpha_n}
(\vert X \vert- \vert a \vert )^{r+\nu} \mathbf{1}_{ \{ \vert X
\vert> \rho_n +c \}}  \Bigr) \nonumber
\\[-8pt]
\\[-8pt]
& \geq&n^{\fracb{r+\nu}{d}} \mathbb{E} \bigl (  (\vert X \vert- \rho
_n  )^{r+\nu} \mathbf{1}_{ \{ \vert X \vert> \rho_n + c \}}
 \bigr)  \nonumber\\
& \geq& c^{r+\nu} n^{\fracb{r+\nu}{d}}\mathbb{P}  ( \vert X \vert>
\rho_n + c  ) .
\nonumber
\end{eqnarray}
Taking the supremum over $c > 0$ and using that $X$ has an $(r,r+\nu
)$-distribution, we complete the proof.

(b) is proved like (a). Inequality (\ref{result_inter})
has the following counterpart: For every $u>1$,
\[
\mathbb{E} \vert X - \widehat{X}^{\alpha_n} \vert^{r+\nu} \geq \mathbb
{E}  \bigl(  (\vert X \vert- \rho_n  )^{r+\nu} \mathbf{1}_{ \{
\vert X \vert> u \rho_n \}}  \bigr) \geq\mathbb{E} \bigl (\vert X
\vert^{r+\nu}  (1- 1/u  )^{r+\nu} \mathbf{1}_{ \{ \vert X
\vert> u \rho_n \}}  \bigr).
\]
Inequality (\ref{eqthm1.6.1}) follows from
\[
n^{\fracb{r+\nu}{d}} \mathbb{E} \vert X - \widehat{X}^{\alpha_n} \vert
^{r+\nu} \geq \sup_{u>1} \bigl [  (1- 1/u  )^{r+\nu} n^{\fracb
{r+\nu}{d}} \mathbb{E}  \bigl(\vert X \vert^{r+\nu} \mathbf{1}_{ \{
\vert X \vert> u \rho_n \}}  \bigr)  \bigr].
\]
\upqed
\end{pf*}

%-----------------------END PROOF------------------
%s4.1.2 ###
\subsubsection{Application to distributions with polynomial or
hyper-exponential tails}
The next proposition is the counterpart of Proposition \ref
{cor_princip_limsup} devoted to the asymptotic lower bound.
\begin{prop} \label{cor_princip_liminf} Let $r>0$ and let $X \in
L^{r+}(\mathbb{P})$ be an $\mathbb{R}^d$-valued random variable having
an unbounded support.
\begin{longlist}[(b)]
\item[(a)] \textup{Polynomial tail.} Set
%
%e4.13 ###
\begin{equation}
\zeta_{\star} = \inf \Bigl\{ \zeta>0,   \forall\nu \in[0,\nu^{\star
}_{X}),   \liminf_{\xi\rightarrow+\infty} \xi^{\zeta-r -\nu} \bar
{F}_{r+\nu}(\xi) >0  \Bigr\}  \in[r+\nu^{\star}_{X}, +\infty].
\end{equation}
Then
%
%e4.14 ###
\begin{equation}\label{3.14}
\liminf_n \frac{\log\rho_n}{\log n} \geq\frac{1}{\zeta_{\star}-r-\nu
^{\star}_{X}} \frac{r+\nu^{\star}_{X}}{d}.
\end{equation}

\item[(b)] \textup{Hyper-exponential tail.} Set
%
%e4.15 ###
\begin{equation} \label{def_expon_tail}
\theta_{\star} = \inf \Bigl\{\theta>0,   \liminf_{\xi\rightarrow
+\infty} \mathrm{e}^{\theta\xi^{\kappa}} \mathbb{P}(\vert X \vert>\xi) >0
\Bigr\}  \in[ 0,+\infty].
\end{equation}
Then, $\theta^\star\le\theta_\star$ and
%
%e4.16 ###
\begin{equation}
\liminf_{n} \frac{\rho_n}{ (\log n  )^{1/ \kappa}} \geq
\biggl(\frac{r + \nu^{\star}_{X}}{d \theta_{\star}}  \biggr)^{1/ \kappa}.
\end{equation}
\end{longlist}
\end{prop}

\begin{pf}
(a) Let $\zeta \in(0, \zeta_{\star
})$. For every $\nu \in(0,\nu^{\star}_{X})$, there exists a positive
real constant~$C_{\nu}$  such that $ \bar{F}_{r+\nu}(\xi) \geq C_{\nu}
\xi^{-\zeta+r+\nu}$ for large enough $\xi$. Setting $\psi_{r,\nu}(y) =
\frac{y}{\zeta-r-\nu}$ yields  $\psi_{r,\nu}(-\log\bar{F}_{r+\nu}(\xi))
\leq\log\xi+\mathrm{o}(\log\xi)$. It follows from Theorem~\ref
{thm_princip_liminf}(a) that
\[
\liminf_{n} \frac{\log\rho_n}{\log n} \geq\frac{1}{\zeta-r-\nu} \frac
{r+\nu}{d}.
\]
Letting $\nu$ and $\zeta$ go to $\nu^{\star}_{X}$ and $\zeta_{\star}$
yields the announced result.

(b) Let $\theta \in (\theta_{\star},+\infty)$. Then, there exists
a positive real constant $C$ such that $ \bar{F}(\xi) \geq C \mathrm{e}^{-\theta
\xi^{\kappa}}\!$  for large enough\vadjust{\goodbreak} $x$. Therefore $ -\log\bar{F}(\xi)
\leq\theta\xi^{\kappa} (1 - \xi^{-\kappa} \log(C) ) $ so that, by
setting $\psi_{\theta}(y) = (y/\theta)^{1/ \kappa}$, we have
\[
\psi_{\theta}(-\log\bar{F}(\xi)) \leq\xi+ \mathrm{o}(\xi).
\]
It follows from Theorem~\ref{thm_princip_liminf}(b) that
\[
\liminf_{n} \frac{\rho_n}{ (\log n  )^{1/ \kappa}} \geq
\biggl(\frac{r + \nu^{\star}_{X}}{\theta d}  \biggr)^{1/ \kappa}.
\]
Letting $\theta$ go to $\theta_{\star}$ completes the proof. Finally,
the inequality between $\theta_\star$ and $\theta^\star$ is an easy
consequence of the fact that $\bar F_r(\xi) \ge\xi^r\bar F(\xi)$.
\end{pf}

Now we give explicit bounds and rates for several families of
distribution tails (which include most usual distributions). To do so,
we combine asymptotic upper bound results from Section~\ref
{Upperexplicit} with asymptotic lower bound results obtained in this
section. The results below are fully explicit in that we make no \textit{a
priori assumptions} on $\nu^{\star}_{_X}$.

\begin{cor} \label{cor_explic_dens_liminf} Let $r>0$ and let $X \in
L^{r+}(\mathbb{P})$ be an $\mathbb{R}^d$-valued random variable, with
probability density $f$, having an unbounded convex support.

\begin{longlist}[(b)]
\item[(a)]\textup{Polynomial tail.} If there exists $c'\ge c>r+d$
such that
\[
0< \liminf_{|x|\to+\infty} |x|^{c'} f(x)\quad\mbox{and}\quad \limsup
_{|x|\to+\infty} |x|^{c} f(x) <+\infty,
\]
then $f$ satisfies  (\ref{LocalControl}),
%
%e4.17 ###
\begin{eqnarray} \label{gen_density_pol_liminf}
&&d \biggl(1-\frac{d+r}{c'} \biggr)-(r+d) \biggl(1-\frac{c}{c'} \biggr)\le\nu
^{\star}_{X} \le d \biggl(1-\frac{d+r}{c'} \biggr),\nonumber
\\[-8pt]
\\[-8pt]
  && \quad  c-d\le \zeta
^{\star},  \qquad  \zeta_\star\le c'-d
\nonumber
\end{eqnarray}
and
\[
\frac{1}{c'-r-d} \biggl(1+\frac rd \biggr)\le \liminf_n \frac{\log\rho
_n}{\log n} \le \limsup_n \frac{\log\rho_n}{\log n}\le \frac
{1}{c-r-d}  \biggl(1+\frac rd \biggr).
\]
Finally, if $c=c'$, then
%
%e4.18 ###
\begin{equation}\label{eq_liminf_limsup_pol}
\nu^{\star}_{X} = d \biggl(1-\frac{d+r}{c'} \biggr), \qquad   \zeta_{\star} =
\zeta^{\star} = c-d\quad\mbox{and} \quad
\lim_{n} \frac{\log\rho_n}{\log n} = \frac{1}{c-r-d}  \biggl(1+\frac
rd \biggr).
\end{equation}
\item[(b)] \textup{Hyper-exponential tail.} If there exists $\kappa>0$ such that
%
%e4.19 ###
\begin{equation} \label{gen_density_exp_liminf}
\lim_{|x|\to+\infty}\frac{\log f(x)}{|x|^{\kappa}} =- \vartheta \in
(-\infty,0),
\end{equation}
then
\[
\nu^{\star}_{X} = d\quad\mbox{and}\quad \theta_{\star} =\theta
^{\star} = \vartheta\vadjust{\goodbreak}
\]
so that
%
%e4.20 ###
\begin{equation} \label{eq_liminf_limsup_exp}
\frac{1}{ \vartheta^{1/ \kappa}}  \biggl(1+\frac{r }{d}  \biggr)^{1/ \kappa
} \leq\liminf_{n} \frac{\rho_n}{ (\log n  )^{1/ \kappa}} \leq
\limsup_{n} \frac{\rho_n}{ (\log n  )^{1/ \kappa}} \leq\frac{2}{
\vartheta^{1/ \kappa}}  \biggl(1+\frac{r}{d}  \biggr)^{1/ \kappa}.\quad
\end{equation}
When $d=1$, $ r \geq1$,
then the following sharp rate holds
%
%e4.21 ###
\begin{equation} \label{eq_liminf_limsup_expdim1}
\lim_{n} \frac{\rho_n}{ (\log n  )^{1/ \kappa}} =  \biggl(\frac
{r+1}{\vartheta}  \biggr)^{1/ \kappa}.
\end{equation}
\end{longlist}
\end{cor}

\begin{remark*} When $d=1$, a one-sided result follows by
considering ``$x\to+ \infty$'' instead of ``$|x|\to\infty$''.
\end{remark*}

\begin{pf*}{Proof of Corollary~\ref{cor_explic_dens_liminf}}
(a) First we need to check that $f$ satisfies the control
criterion~(\ref{LocalControl}) from Proposition~\ref{Critere2}: Let
$A$, $A'$ and $B$ be such that $A'|x|^{-c'}\le f(x)\le A |x|^{-c}$ for
every $x \in\mathbb{R}^d$,    $|x|\ge B$. Then, if $\eta \in(0,1)$,
one checks that the criterion is satisfied with $M= \frac{B}{1-\eta}$,
$K= \frac{A'}{A^{{c'}/{c}}}(1+\eta)^{-c'}$ and $\varepsilon=\frac
{c'-c}{c}\ge0$.

Using that $A'|x|^{-c'}\le f(x)$ and that $f$ is a probability density
(so that $f^a$ is locally integrable if $a \in(0,1]$) yields by
checking (\ref{necess_rate_opt}) the upper bound for $\nu^\star_{_X}$.
Checking now the integral criterion (\ref{Hypoth2})
yields the lower bound.

The lower bound for $\zeta^\star$ is established in Corollary~\ref
{cor_gen_density}. The upper-bound is obtained by similar computations
that show that, if $\zeta>c'-d$, then for $\xi$ large enough, $\xi
^{\zeta-r}\bar F_r(\xi) \ge Ad V_d\xi^{\zeta-(c'-d)}$ for some real
constant $A>0$. This shows that $\zeta^{\star} \le c'-d$.
The bounds for~$\zeta_\star$ are obtained by similar computations.

As concerns the lower bound for the radius, one concludes by plugging
all these estimates into (\ref{3.14}). Combining this with
Corollary~\ref{cor_gen_density}(a) completes this part of the proof.

(b) First we need to check that $f$ satisfies the control
criterion (\ref{LocalControl}). We know from assumption (\ref
{gen_density_exp_liminf}) that for every $ \bar\eta \in(0,\vartheta
)$, there exists $B_{\bar\eta}>0$ such that $\mathrm{e}^{-(\vartheta+\bar\eta
)|x|^{\kappa}}\le f(x)\le \mathrm{e}^{(-\vartheta+\bar\eta)|x|^{\kappa}}$, as
soon as $|x|\ge B_{\bar\eta}$. Then, one shows that the criterion is
satisfied with $M=\frac{B_{\bar\eta}}{1-\eta}$, $K=1$, $\varepsilon
=\frac{\vartheta+\bar\eta}{\vartheta-\bar\eta}(1+\eta)^{\kappa}-1$.
Then, one checks that $\nu^{\star}_{_X}\ge d-(r+d)\frac{\varepsilon
}{1+\varepsilon}$ since $\int_{\{|x|\ge B \}} \exp{(-\mu|x|^{\kappa
})}\,\mathrm{d}x<+\infty$ for every $B$, $\mu>0$. Letting $\eta$ and $\bar\eta
\to0$ yields $\nu^{\star}_{_X}=d$.

To compute $\theta_{\star}$, one first notes that, as soon as $\xi\ge
B_{\bar\eta}$,
\begin{eqnarray*}
\mathbb{P}(|X|>\xi) &\ge& d V_d \int_{\{u>\xi\}} \mathrm{e}^{-(\vartheta+\bar
\eta)u^{\kappa}}u^{d-1}\,\mathrm{d}u\\
&=& \mathrm{O}\bigl(\mathrm{e}^{-(\vartheta+\bar\eta)\xi^{\kappa}}\xi^{d-\kappa}\bigr),
\end{eqnarray*}
where the equality follows by a standard argument based on an
integration by parts and a comparison theorem for integrals. As a
consequence $\theta_{\star}\le\vartheta+\bar\eta$, which finally
implies $\theta_{\star}\le\vartheta$. Combining this
with Corollary~\ref{cor_gen_density}(b) and Proposition~\ref
{cor_princip_liminf}(b) yields $\theta_{\star}=\theta^{\star
}=\vartheta$.
\end{pf*}

\begin{pf*}{Proof of Theorem~\ref{1.2}} Claim (a) follows
from the former Corollary~\ref{cor_explic_dens_liminf}(a) once it is
noted that for every $\varepsilon \in(0,c)$, $\liminf_{|x|\to\infty
}|x|^{c+\varepsilon}f(x)>0$ and $\limsup_{|x|\to\infty}
|x|^{c-\varepsilon}f(x)<+\infty$. Claim~(b) directly follows from
(b) in the above corollary.
\end{pf*}

%---------------------RANDOM QUANTIZATION APPROACH-----------------

%s4.2 ###
\subsection{An alternative approach based on random quantization}\label
{loweriid}
Random quantization is another tool to compute the lower estimate of
the maximal radius sequence of a random vector $X$ with distribution
$P$. It makes a connection between $\rho_n$ and the maximum of an
i.i.d. sequence of random variables with distribution $P$.
\begin{theo} \label{theoreme1} Let $r>0$ and let $ X \in L^{r+}(\mathbb
{P}) $ be a random variable taking values in $\mathbb{R}^d$ with an
absolutely continuous distribution $P$.
Assume $(\alpha_n)_{n \geq1}$ is a sequence of $L^r(P)$-optimal
$n$-quantizers. Let $( X_k)_{k \geq1}$ be an   i.i.d.  sequence of
$\mathbb{R}^d$-valued copies of $X$. For every $\nu \in[0,\nu^{\star
}_{X})$ such that $r+\nu\ge1$, there exists a real constant $C_{r,\nu}
\in(0,\infty)$ such that\looseness=1
%
%e4.22 ###
\begin{equation} \label{first_low_estim}
\liminf_{n}  \Bigl( \rho_n - \mathbb{E}  \Bigl( \max_{k \leq[n^{(r+\nu
)/d}]} \vert X_k \vert \Bigr)  \Bigr) \geq - C_{r,\nu}.
\end{equation}\looseness=0
\end{theo}
\begin{pf}  Let $\nu \in[0,\nu^{\star}_{X})$ and
set $\widehat{X}_{k}^{\alpha_n} = \sum_{a \in\alpha_n} a \mathbf{1}_{\{ X_k \in C_a (\alpha_n) \}}$. We have, for integer $m \ge1$,
\begin{eqnarray*}
\rho_n & \geq & {\max_{k \leq m}} \vert\widehat{X}_{k}^{\alpha_n}
\vert \\
& \geq& \sum_{k=1}^{m} {\max_{l \leq m}} \vert\widehat{X}_{l}^{\alpha
_n} \vert\mathbf{1}_{\{ \vert X_k \vert> \max\{ \vert X_i \vert,
i \in\{1,\ldots ,m \}, i\neq k \} \}}\\
& \geq& \sum_{k=1}^{m} \vert\widehat{X}_{k}^{\alpha_n} \vert\mathbf{1}_{\{ \vert X_k \vert> \max_{i \not= k} \vert X_i \vert\}} \\
& \geq& \sum_{k=1}^{m}  (\vert X_k \vert- \vert X_k - \widehat
{X}_{k}^{\alpha_n} \vert ) \mathbf{1}_{\{ \vert X_k \vert> {\max
_{i \not= k}} \vert X_i \vert\}}.
\end{eqnarray*}
Taking the expectation of both sides of the previous inequality yields
\[
\rho_n \geq \mathbb{E} {\max_{k \leq m}} \vert X_k \vert- \sum
_{k=1}^{m} \mathbb{E}  \bigl(\vert X_k - \widehat{X}_{k}^{\alpha_n}
\vert\mathbf{1}_{\{ \vert X_k \vert> {\max_{i \not= k}} \vert X_i
\vert\}}  \bigr).
\]
Furthermore, $ \forall k \geq1$, $ \vert X_k - \widehat{X}_k^{\alpha
_n} \vert\mathbf{1}_{\{ \vert X_k \vert> \max_{i \not= k} \vert X_i
\vert\}}$ and $\vert X_1 - \widehat{X}_1^{\alpha_n} \vert\mathbf{1}_{\{ \vert X_1 \vert> \max_{i \not= 1} \vert X_i \vert\}}$ have the
same distribution. Hence,
\begin{eqnarray*}
\rho_n & \geq & \mathbb{E} {\max_{k \leq m}} \vert X_k \vert- m \mathbb
{E}  \bigl(\vert X_1 - \widehat{X}_{1}^{\alpha_n} \vert\mathbf{1}_{\{ \vert X_1 \vert> {\max_{i \not= 1}} \vert X_i \vert\}} \bigr) \\
& \geq& \mathbb{E} {\max_{k \leq m}} \vert X_k \vert- m \Vert X_1 -
\widehat{X}^{\alpha_n}_1 \Vert_{r+\nu} \Bigl (\mathbb{P}  \Bigl(\vert X_1
\vert> {\max_{i \not= 1}}\vert X_i \vert \Bigr)  \Bigr)^{1-1/(r+\nu)}
\end{eqnarray*}
owing to the H\"{o}lder inequality. Since the events $ \{ \vert X_k
\vert> {\max_{i \not= k}} \vert X_i \vert\},  k=1,\ldots ,m $,
are pairwise disjoint with the same probability, we have $\mathbb{P}
 (\vert X_1 \vert> {\max_{i \not= 1}}\vert X_i \vert )\leq\frac{1}{m}.$
Finally,
\[
\rho_n \geq \mathbb{E} {\max_{k \leq m}} \vert X_k \vert- m^{\fracc
{1}{r+\nu}} \Vert X - \widehat{X}^{\alpha_n} \Vert_{r+\nu}.
\]
It follows, by setting $ m = [n^{(r+\nu)/d}] $, that
\[
\liminf_{n} \Bigl ( \rho_n - \mathbb{E}  \Bigl( {\max_{k \leq[n^{(r+\nu
)/d}]}} \vert X_k \vert \Bigr)  \Bigr) \geq- \limsup_{n} n^{
{1}/{d}}\Vert X - \widehat{X}^{\alpha_n} \Vert_{r+\nu}.\vspace*{-2pt}
\]
The upper limit on the right-hand side is finite since $X$ has an
$(r,r+\nu)$-distribution.\vspace*{-2pt}
\end{pf}

\begin{exam} [(Exponential distribution)]  Let $r>0$ and let $X$ be
an exponentially distributed random variable with parameter $\lambda
>0$. If $(\alpha_n)_{n \geq1}$ is an $L^r$-optimal sequence of
$n$-quantizers for~$X,$ then Theorem \ref{theoreme1} implies
%
%e4.23 ###
\begin{equation} \label{eq_exp_comment}
\liminf_{n} \frac{\rho_n}{ \log n} \geq\frac{r+1}{\lambda},\vspace*{-2pt}
\end{equation}
which corresponds to the sharp rates given by (\ref{asymp_exp})
and (\ref{eq_liminf_limsup_pol}), respectively.\vspace*{-2pt}
\end{exam}

Indeed, let $\nu \in(0,\nu^{\star}_{X})$ and let $(X_i)_{i=1,\ldots
,[n^{r+\nu}]}$, be an i.i.d. exponentially distributed sequence of
random variables with parameter $\lambda$. We have for every $u \geq0$,
\[
\mathbb{P}\Bigl(\max_{i \leq[n^{r+\nu}]} X_i \geq u\Bigr) = 1 - \mathbb{P}(X
\leq u)^{[n^{r+\nu}]} = 1 - F(u)^{[n^{r+\nu}]},\vspace*{-2pt}
\]
where $F$ is the distribution function of $X$ (we will denote by $f$
its density function).~Then
\begin{eqnarray*}
\mathbb{E} \Bigl ( \max_{i \leq[n^{r+\nu}]} X_i  \Bigr)
& = & \int_0^{+ \infty} \mathbb{P}\Bigl(\max_{i \leq[n^{r+\nu}]} X_i \geq
u\Bigr)\,\mathrm{d}u = \int_0^{+\infty} \bigl(1 - (1-\mathrm{e}^{- \lambda u})^{[n^{r+\nu}]}\bigr)\,\mathrm{d}u \\[-2pt]
& = & \int_0^{+\infty} \bigl (1 + F(u) + \cdots+ F(u)^{[n^{r+\nu}]-1}
 \bigr) \frac{f(u)}{\lambda}\,\mathrm{d}u \\[-2pt]
& = & \frac{1}{\lambda} \biggl(1 + \frac{1}{2} + \cdots+ \frac{1}{[n^{r+\nu
}]}\biggr) \\[-2pt]
& \geq & \frac{1}{\lambda} \log(1+ [n^{r+\nu}]) \geq \frac{r+\nu
}{\lambda} \log n.\vspace*{-2pt}
\end{eqnarray*}
Consequently, it follows from the super-additivity of the liminf that
for every $\nu \in(0,1)$,
\[
\liminf_{n} \frac{\rho_n}{\log n}  \geq \liminf_{n} \frac{ \rho_n -
\mathbb{E}  ( \max_{i \leq[n^{r+\nu}] } X_i  )}{\log n} + \liminf
_{n} \frac{\mathbb{E}  ( \max_{i \leq[n^{r+\nu}]}X_i  )}{ \log
n}  \geq \frac{r+ \nu}{\lambda}.\vspace*{-2pt}
\]
The result follows by letting $\nu$ go to $\nu^{\star}_{X}=1$.

In fact, one may easily extend this example to a more general
framework, although, overall, the connection made in Theorem~\ref
{theoreme1} seems less straightforward in terms of deriving explicit
asymptotic lower bounds than the former approach based on more
geometric arguments.\vspace*{-2pt}

\begin{exam} [(Radial distribution with exponential tails)] Let $X$
be an $\mathbb{R}^d$-valued random vector with an unbounded support
having an absolutely continuous distribution with a radial probability
density $f(x)=g(|x|_{S})$ with respect to an Euclidean norm $|\cdot
|_{S}$ so that $\bar F(\xi)= K_{d,S} \int_{\xi}^{+\infty}
u^{d-1}g(u)\,\mathrm{d}u$, $\xi>0$, with $ K_{d,S}= d V_d
(\operatorname{det}(S))^{-
1/2}>0$.\vadjust{\goodbreak} Assume that $\bar F(\xi)\ge cf(\xi)$ for $\xi\ge A>0$ for some
real constant $c>0$. Then
\[
\liminf_{n} \frac{\rho_n}{\log n} \ge c(r+\nu_{_X}^{\star}).\vspace*{-2pt}
\]
\end{exam}

\begin{exam} [(Pareto distribution)] Let $X$ be a random variable
having a Pareto distribution with index $\gamma>0 $. If $(\alpha_n)_{n
\geq1}$ is an asymptotically $L^r$-optimal sequence of $n$-quantizers
for $X$, $r \in(0, \gamma)$, then Theorem \ref{theoreme1} yields
\[
\liminf_{n} \frac{\log\rho_n}{ \log n} \geq\frac{r+1}{\gamma+1},\vspace*{-2pt}
\]
which is not the sharp rate given by (\ref{asymp_Par}).\vspace*{-2pt}
\end{exam}

Notice that if $\gamma>r$, then $X \in L^{r+\eta}( \mathbb{P})$ for $
\eta \in (0,\gamma- r  )$. Now, to prove this result, let $\nu
 \in(0, \nu^{\star}_{X})$ and let $(X_i)_{i\ge1}$ be an i.i.d.
sequence of Pareto-distributed random variables (with index $\gamma$).
We have
\[
\forall m \geq1, \forall u \geq1  \qquad \mathbb{P} \Bigl(\max_{i \leq m}
X_i \leq u \Bigr) = (1-u^{- \gamma})^{m}.\vspace*{-2pt}
\]
Then, the density function of $\max_{1 \leq i \leq m} X_i$ is $ m
\gamma u^{-(\gamma+1)}(1 - u^{-\gamma})^{m-1}.$ Hence,
\begin{eqnarray*}
\mathbb{E}  \Bigl( \max_{1 \leq i \leq m } X_i  \Bigr) & = & m \gamma \int
_1^{+ \infty} x^{- \gamma} (1-x^{-\gamma})^{m-1}\,\mathrm{d}x
=m   B\biggl(1- \frac{1}{ \gamma},m\biggr) \\[-2pt]
& = & \frac{\Gamma(1-  {1}/{ \gamma})\Gamma(m+1)}{ \Gamma(m+1-
{1}/{\gamma})}
\sim\Gamma\biggl(1- \frac{1}{ \gamma}\biggr) m^{ {1}/{\gamma}} \qquad\mbox{as }   m \rightarrow+\infty,\vspace*{-2pt}
\end{eqnarray*}
where we used Stirling's formula for the last statement ($B(\cdot,\cdot)$
denotes the beta function of the first kind). We finally set $m =
[n^{r+\nu}] $ to get
\[
\mathbb{E} \Bigl ( \max_{1 \leq i \leq[n^{r+\nu}] } X_i  \Bigr) \sim
\Gamma\biggl(1- \frac{1}{ \gamma}\biggr) n^{\fracb{r+\nu}{\gamma}}.\vspace*{-2pt}
\]
It follows from (\ref{first_low_estim}) that for every $\varepsilon
\in(0,1)$, $ \rho_n - (1-\varepsilon) \Gamma(1- \frac{1}{ \gamma})
n^{\fracb{r+\nu}{\gamma}} \geq-( \mathrm{C}_{r,\nu}+\varepsilon)$.
Dividing both sides of the inequality by $n^{\fracb{r+\nu}{\gamma}}$ and
taking the logarithm yields
\[
\log\rho_n - \frac{r+\nu}{\gamma} \log n \geq\log \biggl( (1-\varepsilon
 ) \Gamma\biggl(1- \frac{1}{ \gamma}\biggr) - (\varepsilon+ \mathrm{C}_{r,\nu})
n^{-\fracb{r+\nu}{\gamma}}  \biggr).\vspace*{-2pt}
\]
Consequently ${\liminf}_{n \rightarrow+\infty} \frac{\log\rho
_n}{\log n} \geq\frac{r+\nu}{\gamma}$ for every $\nu \in(0,\nu^{\star
}_{X})$. One concludes by letting $\nu$ go to $\nu^{\star}_{X}= \frac
{\gamma-r}{\gamma+1}$.\vspace*{-2pt}

\begin{comment*} Let $\phi$ be the inverse (if any)
function of $-\log\bar{F}$. Notice that in both examples above we have
%
%e4.24 ###
\begin{equation} \label{eq_par_comment}
\lim_n \frac{ \mathbb{E}  ({\max_{k \leq[n^{r+\nu^{\star}_{X}}]}}
\vert X_k \vert )}{\phi((r+\nu^{\star}_{X}) \log n)}=1,\vspace*{-2pt}\vadjust{\goodbreak}
\end{equation}
which, for distributions with hyper-exponential tails, leads to the
asymptotic lower bound~(\ref{eq_thm_gen_liminf1}) for the sequence
$(\rho_n)_{n\geq1}$. As concerns Pareto distribution, using the
approximation~(\ref{eq_par_comment}) to compute the asymptotic lower
estimate of the maximal radius sequence induces the loss of the
``$-r$'' term in the exact asymptotics. To recover this remaining term
we have simply to consider the inverse function of $-\log\bar{F}_{r+\nu
^{\star}_{X}}$ (as done in the previous section) instead of $-\log\bar
{F}$, and the random quantization approach clearly does not allow us to
do so.
\end{comment*}

%s4.2.1 ###
\subsubsection{A conjecture about the sharp rate}\label{sec4.2.1}
The previous results related to distributions with hyper-exponential
tails strongly suggest the following conjecture: Suppose $X$ is a
distribution with hyper-exponential tail in the sense of statement (\ref
{gen_density_exp_liminf}). Then, for every $r>0$ and $d \geq1$,
\[
\lim_{n} \frac{\rho_n}{ (\log n  )^{1/ \kappa}} =  \biggl(\frac
{r+d}{d \theta^{\star}}  \biggr)^{1/ \kappa}.
\]
This conjecture is proved for $d=1$ and $r\geq1$. To be satisfied for
higher dimensions we need to prove that the geometric statement (\ref
{equaconjecture}) of Lemma \ref{lem1_princip_limsup} holds true with
$c_{r,d}=1$ for every $r>0$, $d \geq1$.

\begin{appendix}\label{appm}
\section*{Appendix}

\mbox{}\vspace*{-10pt}

\setcounter{equation}{0}
$\rhd$ \textit{Exponential distribution.} $\rho_n = \frac
{r+1}{\lambda} \log n+ \frac{C_{r}}{\lambda} + \mathrm{O} (\frac{1}{n}
 ),$ we use the following result (see \cite{ForPag}): If $X$ is
exponentially distributed with parameter $\lambda>0,$ then, for any $n
\geq1$, the $L^r$-optimal quantizer $\alpha_n = (\alpha_{n,1}, \ldots ,
\alpha_{n,n})$ is unique and given by
%
%e4.1 ###
\begin{equation}
\alpha_{n,k} = \frac{1}{\lambda}\Biggl (\frac{a_n}{2} + \sum
_{i=n+1-k}^{n-1} a_i  \Biggr), \qquad 1 \leq k \leq n,
\end{equation}
where $(a_k)_{k \geq1}$ is an $\mathbb{R}_{+}$-valued sequence
recursively defined by the following implicit equation: $a_0:= +\infty,
 \phi_r(-a_{k+1}):= \phi_r(a_k), k \geq0 $, with $\phi_r(x):= {\int
_0^{x/2}} \vert u \vert^{r-1} \operatorname{sign}(u) \mathrm{e}^{-u}\,\mathrm{d}u $ (convention: $
0^0 = 1$). Furthermore, the sequence $(a_k)_{k \geq1}$ decreases to
zero and for every $k \geq1$, $a_k = \frac{r+1}{k} (1+ \frac
{c_r}{k} + \mathrm{O}(\frac{1}{k^2})  ) $ for some positive real
constant $c_r$. Then it follows that $\lambda\rho_n = \frac{a_n}{2} +
\sum_{i=1}^{n-1} a_i $ so that
\[
\lambda\rho_n = \frac{a_n}{2} + (r+1) \sum_{i=1}^{n-1} \frac{1}{i} +
c_r \sum_{i=1}^{n-1} \frac{1}{i^2} + \sum_{i=1}^{n-1} \mathrm{O}(1/i^3) =
(r+1) \log n + C_r + \mathrm{O} \biggl(\frac{1}{n}  \biggr).
\]

$\rhd$ \textit{Pareto distribution.} The proof is similar after
noting that $ \rho_n = \frac{1}{1+a_n} \prod_{i=1}^{n-1} (1+a_i),$
where $(a_n)_{n \geq1}$ is an $\mathbb{R}_{+}$-valued sequence,
decreasing to zero and satisfying, for every $n \geq1$, $ a_n = \frac
{r+1}{(\gamma- r)n}  ( 1 + c_r /n + \mathrm{O}(1/n^2)  )$ for
some\vadjust{\goodbreak}
real constant $ c_r $. Hence,  if $i_0:=\max\{i \mid  |a_i|\ge1\}$,
\[
\log(\rho_n) = -\log(1+a_n) +C_{i_0}+ \sum_{i=i_0+1}^{n-1}  \biggl( a_i -
\frac{a_i^2}{2} + \mathrm{O}(a_i^3)  \biggr) = \frac{r+1}{\gamma-r} \log n +
C_r + \mathrm{O}  \biggl(\frac{1}{n}  \biggr),
\]
where we used that $\sum_{i=1}^{\infty} a_i^2 < \infty$ and $\sum
_{i=1}^{\infty} \mathrm{O}(a_i^3) < \infty$.
\end{appendix}

% imsref loaded by smiklovaite, 2011-03-18 09:35:29
% imsref loaded by smiklovaite, 2011-03-18 10:01:58
%

\printhistory

\end{document}